\documentclass[11pt]{article}

\usepackage[top=50mm, bottom=50mm, left=50mm, right=50mm]{geometry}
\usepackage{lineno}                 % Line numbers
\usepackage{amssymb}                % Math
\usepackage{amsmath}                % Math
\usepackage{amsthm}                 % Math
\usepackage{mathrsfs}               % Math, \mathscr{}
\usepackage{epsfig}                 % eps figures
\usepackage{graphicx}               % Figures
\usepackage{graphics}               % ?
\usepackage{float}                  % 
\usepackage{multirow}               % 
\usepackage{color}                  %
\usepackage{fullpage}               %
\usepackage[normalem]{ulem}         % 
\usepackage{makeidx}                %
\usepackage{xspace}                 %
\usepackage{wrapfig}                %

\usepackage{url}                    % simple URL typesetting
\usepackage{booktabs}               % professional-quality tables
\usepackage{amsfonts}               % blackboard math symbols
\usepackage{nicefrac}               % compact symbols for 1/2, etc.
\usepackage{microtype}              % microtypography

\usepackage{caption}                %
\usepackage{enumitem}               %
\usepackage[export]{adjustbox}      %
\usepackage{comment}                %
\usepackage{xstring}                %
\usepackage{cite}                   %
\usepackage{amsfonts}               %
\usepackage{subcaption}             %
\usepackage{tikz}                   % 
\usepackage{pgfplots}               %
\usepackage{pgfplotstable}          %

\usepackage{multicol}               % Multi-columns

\usepackage[affil-it]{authblk}                % author and affiliations 

\usepackage[ruled, vlined, linesnumbered]{algorithm2e}
\usepackage{mathtools}

\makeindex

\pgfplotsset{compat=1.17}

\newtheorem{problem}{Problem}
\newtheorem{dproblem}{Discrete Problem}
\newtheorem{remark}{Remark}

\def\bepsilon{\hbox{\boldmath$\epsilon$}}
\def\bnabla{\hbox{\boldmath$\nabla$}}

\def\bdisp{\hbox{$\mathbf{u}$}}

\def\testbdisp{\hbox{$\delta\mathbf{u}$}}

\def\strain{\hbox{$\bepsilon$}}

\newcommand{\lame}[1]{%
    \IfEqCase{#1}{%
        {1}{\lambda}%
        {2}{\mu}%
    }[\PackageError{lame}{Undefined option to lame: #1}{}]%
}% Lame constants

\def\pf{\hbox{$\varphi$}}
\def\testpf{\hbox{$\delta\varphi$}}
\def\l{\hbox{$l$}}
\def\gc{\hbox{$G_c$}}

\def\dom{\hbox{$\Omega$}}
\def\surf{\hbox{$\Gamma$}}

\newcommand{\surfarg}[2]{\hbox{$\surf_{#1}^{\scriptsize{#2}}$}}
\newcommand{\trac}[1]{\hbox{$\mathbf{t}_{p}^{\scriptsize{#1}}$}}

\def\surfcrack{\hbox{$\mathcal{C}$}}

\usetikzlibrary{fadings}
\tikzfading[name=fade out, 
    inner color=transparent!0,
    outer color=transparent!100]
\usetikzlibrary{shapes}
\usetikzlibrary{decorations.pathreplacing}
\usetikzlibrary{arrows}
\pgfplotsset{
        colormap={parula}{
            rgb255=(53,42,135)
            rgb255=(15,92,221)
            rgb255=(18,125,216)
            rgb255=(7,156,207)
            rgb255=(21,177,180)
            rgb255=(89,189,140)
            rgb255=(165,190,107)
            rgb255=(225,185,82)
            rgb255=(252,206,46)
            rgb255=(249,251,14)
        },
    }
    
\providecommand{\keywords}[1]
{
  \small	
  \textbf{\textit{Keywords: }} #1
}

\begin{document}
%\title{PHT spline-based adaptive isogeometric analysis of phase-field fractures with a novel monolithic solver}

%\title{A robust monolithic solver for phase-field brittle fracture using fracture energy arc-length method and adaptive under-relaxation}

%\title{A robust monolithic implementation of phase-field fracture: fracture energy arc-length method and adaptive under-relaxation}

%\title{ and under-relaxation integrated robust monolithic implementation of phase-field fracture}

%\title{Fracture energy based arc-length method and under-relaxation integrated robust monolithic solver for phase-field fracture}

\title{A robust monolithic solver for phase-field fracture integrated with fracture energy based arc-length method and under-relaxation}

\author[1]{Ritukesh Bharali}
\author[2]{Somdatta Goswami}
%\author[1]{Fredrik Larsson}
%\author[3]{Ralf J\"anicke}
\author[3]{Cosmin Anitescu}
\author[3]{Timon Rabczuk}
\affil[1]{Department of Industrial and Material Science, Chalmers University of Technology}
\affil[2]{Department of Applied Mathematics, Brown University, U.S.A}
%\affil[3]{Institute of Applied Mechanics, Technische Universit\"at Braunschweig}
\affil[3]{Institute of Structural Mechanics, Bauhaus-Universit\"at Weimar}

\date{}

\maketitle

\section*{Highlights}
\begin{itemize}
    \item Robust monolithic solver with adaptive under-relaxation and arc-length method.
    \item Snap-back behaviour is captured with a phase-field fracture energy-based dissipation constraint.
    \item PHT-splines within IGA framework is utilized for adaptive mesh refinement.
\end{itemize}

\bigskip

\begin{abstract}

%\noindent The non-convex nature of the phase-field fracture energy functional renders the monolithic solution techniques non-robust. It is observed that the monolithic solution techniques diverge at the onset of the localisation of the fracture zone. In this manuscript, a novel monolithic solution technique is proposed for the phase-field fracture model, that enables a numerical model to go beyond the onset of fracture. The solution technique utilizes an adaptive under-relaxation method for circumventing the divergence of the solver, and also incorporates a fracture energy dissipation-based arc-length for adaptive stepping of the quasi-static load increments. In order to improve the efficiency of the model, adaptive mesh refinement is considered using the PHT splines. Numerical experiments on benchmark problems are presented to demonstrate the robustness of the proposed solution technique.
\noindent The phase-field fracture free-energy functional is non-convex  with respect to the displacement and the phase field. This results in a poor performance of the conventional monolithic solvers like the Newton-Raphson method. In order to circumvent this issue, researchers opt for the alternate minimization (staggered) solvers. Staggered solvers are robust for the phase-field based fracture simulations as the displacement and the phase-field sub-problems are convex in nature. Nevertheless, the staggered solver requires very large number of iterations (of the order of thousands) to converge. In this work, a robust monolithic solver is presented for the phase-field fracture problem. The solver adopts a fracture energy-based arc-length method and an adaptive under-relaxation scheme. The arc-length method enables the simulation to overcome critical points (snap-back, snap-through instabilities) during the loading of a specimen. The use of an under-relaxation scheme stabilizes the solver by preventing the divergence due to an ill-behaving stiffness matrix. The efficiency of the proposed solver is further amplified with an adaptive mesh refinement scheme based on PHT-splines within the framework of isogeometric analysis. The numerical examples presented in the manuscript demonstrates the efficacy of the solver. All the codes and data-sets accompanying this work will be made available on GitHub (\url{https://github.com/rbharali/IGAFrac}).
% \noindent The non-convex nature of the free-energy functional in phase field based fracture analysis often renders the conventional monolithic solvers incapacitated. To avoid convergence issues at the onset of the fracture localisation, researchers often resort to staggered solvers based on algorithmic decoupling, which are typically slow converging. In this work, we propose a robust and efficient monolithic solver that competently resolves the divergence issues of a conventional monolithic solver (with Ne) and is also capable of capturing the discontinuous response, such as the snap-back instability, which is particularly prominent in displacement-controlled fracture, and are commonly disregarded in staggered solution schemes. The proposed monolithic solver exploits an adaptive under-relaxation method for circumventing the divergence of the solver. Alongside, the proposed solver also incorporates a fracture energy dissipation-based arc-length technique for adaptive displacement increments. The efficiency of the proposed solver is amplified by implementing the adaptive scheme using PHT-splines within the framework of IGA. Numerical experiments on benchmark problems are presented along with the force-displacement curves, to demonstrate the robustness and accuracy of the novel monolithic solver. All the codes and data-sets accompanying this work will be made available on GitHub.
%\SG{mention about large dissipation steps}
\end{abstract}

% The phase-field fracture free energy functional is non-convex  w.r.t the displacement and the phase-field. This renders a poor performance of the conventional monolithic solvers like the Newton-Raphson method. In order to circumvent this issue, researchers opt for the alternate minimization (staggered) solvers. Staggered solvers are robust for the phase-field fracture simulations as the displacement and phase-field sub-problems are convex in nature. Nevertheless, the staggered solver requires very large number of iterations (around thousands) to converge. In this work, a robust monolithic solver is presented for the phase-field fracture problem. The solver adopts a fracture energy-based arc-length method and an adaptive under-relaxation scheme. The arc-length method enables the simulation to overcome critical points (snap-back, snap-through instabilities) during the loading of a specimen. The use of an under-relaxation scheme stabilizes the solver by preventing the divergence due to an ill-behaving stiffness matrix. Through numerical examples presented in this manuscript, it is demonstrated that the proposed arc-length method in conjunction with the adaptive under-relaxation scheme allow large dissipation steps. 

% keywords can be removed
\keywords{phase-field fracture, brittle material, monolithic solver, arc length method, variational damage, IGA, PHT-splines}

\newpage

\section{Introduction}\label{sec:1}

The seminal work of Francfort and Marigo \cite{Francfort1998} led to the emergence of the phase-field based fracture model, as an alternative fracture modelling technique. Therein, the Griffith fracture criterion was cast into a variational setting with certain limitations: no concept of internal length scale and maximum allowable stresses. Later, \cite{Bourdin2000,Bourdin2007} proposed a computationally convenient framework of the Francfort and Marigo model, adopting a scalar auxiliary variable that interpolates between fully fractured and intact material states. In this context, the Ambrosio-Tortorelli regularization of the Mumford-Shah functional \cite{mumford1989optimal} was utilized. Based on the minimization of the global energy function, the phase-field model eliminates the need for tedious tracking of the fracture path and remeshing techniques, frequently observed in the discrete fracture models like XFEM\cite{bordas2011performance}. Furthermore, the phase-field model for fracture has proven its capabilities to handle topologically complex fracture patterns (branching, kinking and merging of cracks)\cite{Bourdin2000}.

Soon after the inception of phase field based fracture model, the concept was cast into a thermodynamically consistent framework in \cite{miehe2010b}, adopting an energy-based fracture driving criterion. This work was later extended towards a generic stress-based fracture driving criterion in \cite{miehe2015449}, ductile fracture with plasticity models in \cite{ambati2015ductile,miehe2015486}, anisotropic fracture \cite{TEICHTMEISTER20171,BLEYER2018213}, hydraulic fracture \cite{WILSON2016264,chukwudozie2019variational}, desiccation cracking \cite{cajuhi2018phase,HEIDER2020112647} in a non-exhaustive list of single-scale brittle fracture applications. In the context of multi-scale modelling, the overlapping domain decomposition techniques were adopted in \cite{patil2018adaptive,nguyen2019multiscale,triantafyllou2020generalized,gerasimov2018non,noii2020adaptive}, while \cite{fantoni2020phase,he2020numerical,Bharali2021} adopted the hierarchical modelling technique in the FE$^2$ sense \cite{feyel2000fe2}.

Despite its popularity in several multi-physics domains, the phase-field model has its own set of computational challenges in its implementation in fracture analysis. They include,
(\textit{1}) a non-convex free-energy functional with respect to the coupled field variables, (\textit{2}) variational inequality due to fracture irreversibility constraint, and (\textit{3}) the need to resolve the smeared fracture zone with an extremely fine mesh. The coupled fields can be solved using either a monolithic solver or a staggered solver. Provided that the non-linear solver converges, the monolithic solution scheme is more efficient and faster than the staggered one. However, the non-convex energy functional generally leads to poor convergence and loss of robustness of the monolithic solver. In order to circumvent this, \cite{Gerasimov2016} proposed a line search technique which included a negative search direction, \cite{Heister2015} proposed convexification of the energy functional based on linear extrapolation of the phase-field for the momentum balance equation. Other methods developed in this context include the arc-length solvers \cite{vignollet2014phase,may2015numerical,singh2016fracture}, modified Newton-Raphson method \cite{wick2017modified}, error-oriented Newton-Raphson method \cite{Wick2017a}, and trust regions methods \cite{kopanivcakova2020recursive}. Nevertheless, the development of a robust monolithic solver still remains an active research area in the phase-field fracture community. As an alternative approach, \cite{Miehe2010} suggested the use of staggered solution scheme, since the energy functional is convex with respect to either of the coupled field, if the other one is held constant.

The second computational challenge associated with the monolithic solver pertains to the variational inequality formulation that arises from the fracture irreversibility constraint. In this context, \cite{Gerasimov2016,GERASIMOV2019990} opted for a simple penalization technique, \cite{Heister2015} adopted the primal-dual active set method, whereas \cite{Wick2017a,wick2017modified} used an Augmented Lagrangian method based on the Moreau-Yoshida indicator function \cite{Wick2017a,wick2017modified}. More recent approaches include the micromorphic approach that transforms the phase-field into a local variable \cite{miehe2016phasemicromorphic,Bharali2022micromorphic}, and the slack variable approach \cite{Bharali2022slack}. Alternatively, a heuristic approach was proposed by \cite{Miehe2010}, replacing the fracture driving energy with its maximum value over the loading history. However, this method is not variationally consistent with the phase-field fracture energy functional \cite{GERASIMOV2019990,DeLorenzis2021}.

Finally, the phase-field model for fracture analysis requires an extremely fine mesh to resolve the smeared region. A simple and straight forward but computationally expensive way would be to use uniformly refined mesh. If the fracture path is known $\it{a\;priori}$, a certain part of the computational domain may also be pre-refined. The latter case is more applicable when it comes to benchmark models from the literature. However, in the scenarios where the fracture path is not known in advance, adaptive mesh refinement techniques is the preferred option. In this context, the elements are marked based on either a critical threshold value defined over the phase field parameter \cite{Heister2015,goswami2020adaptive}, or by local increase in the tensile energy \cite{klinsmann2015assessment}. Other adaptive refinement schemes include the recovery-based error indicator \cite{jansari2019adaptive}, $\it{a\;posteriori}$ error estimation based on the dual-weighted residual method \cite{wick2016goal}, the finite cell method \cite{nagaraja2019phase}, and the dual mesh concept for the two coupled fields with different mesh refinement indicators \cite{goswami2019adaptive}.

To overcome the issues discussed in the current state-of-the-art, we propose a novel monolithic solver, which is based on an adaptive under-relaxation scheme, and is integrated with fracture energy-based arc-length method. Although the under-relaxation strategies result in decreased rate of convergence of the solver, for the phase-field model, it ensures a guarantee to circumvent of divergence issues arising due to erratic behaviour of the jacobian \cite{gerasimov2018non}. The fracture energy-based arc-length method adopted in this work is displacement controlled, which provides the flexibility to take larger displacement steps while accurately capturing the brutal nature of the crack growth. In this work, we have studied displacement controlled fracture, where the load $\it{vs.}$ displacement curve in the post-peak behaviour encounters a discontinuity and the representative point drops onto the lower branch with negative slope, indicating that both load and displacement must decrease to obtain a controlled crack extension. Such observations are often neglected in staggered solvers, but these phenomena are captured accurately using our fast and efficient monolithic solver. Besides capturing the possible snap-back behaviour, the arc-length method also results in an adaptive time-stepping procedure, hence larger energy dissipation is permitted. The adaptive scheme enables a dynamically changing mesh which in turn allows the refinement to remain local at singularities and high gradients. The adaptive $\it{h}$-refinement technique is implemented using polynomial splines over hierarchical T-meshes (PHT-splines). The PHT-splines possess a very efficient local refinement algorithm and they also inherit the properties of adaptivity and locality of T-splines. Moreover, in all the examples, the crack is not initialized but it is allowed to nucleate naturally. The penalization approach is adopted to treat the variational inequality formulation.

This remainder of the article is structured as follows: Section \ref{sec:2} introduces the reader to the phase-field model for fracture analysis, its underlying energy functional and the pertinent Euler-Lagrange equations. Subsequently, in Section \ref{sec:3}, the isogeometric analysis framework and the discrete equations for the phase-field fracture model are introduced. Section \ref{sec:4} presents the main contribution of this work, the robust monolithic solver. The numerical benchmark problems are addressed in Section \ref{sec:5}, followed by concluding remarks in Section \ref{sec:6}.

% ---------------- PHASE-FIELD MODEL FOR FRACTURE -------------------------%
\section{Phase-field model for fracture}\label{sec:2}

\subsection{The energy functional}

\begin{figure*}[ht!]
  \begin{subfigure}[b]{0.49\textwidth}
    \begin{tikzpicture}[scale=0.8]
    \coordinate (K) at (0,0);
    % Potato
    \draw [fill=black!10,line width=1pt] (K) plot [smooth cycle,tension=0.7] %coordinates {(3,1) (5,1.2) (7,1) (8,3) (7,4.5) (5,4.5) (2,4) (1.7,2.5)};
    coordinates {(3,1) (7,1) (8,3) (7,4.475) (5,4.5) (2,4) (1.7,2.5)};
    \node[ ] at (3.25,2.15) {$\dom$};
    % Crack Surface
    \draw[line width=1.5pt,black] (5,0.75) to (5,2.5);
    \node[ ] at (5.35,2.0) {$\surfcrack$};
    % Dirichlet Boundary
    \draw (K) [line width=2.5pt,black] plot [smooth, tension=0.8] coordinates {(3,4.3) (1.7,3.7) (1.7,2.5)};
    %\draw (K) [solid,line width=0.5pt] (1.1,4) - - (1.55,3.5);
    \node[ ] at (1.15,4) {$\surfarg{D}{u}$};
    % Neumann Boundary
    %\draw (K) [line width=2.5pt,black] plot [smooth, tension=0.7] coordinates {
    %(7.95,3.4) (7.05,4.5) (5.05,4.5)};
    %\draw (K) [solid,line width=0.5pt] (6.35,4.5) - - (5.75,4.95);
    \node[ ] at (6.75,5.05) {$\surfarg{N}{u}$};
    \end{tikzpicture}
    \caption{discrete crack}
    \label{fig:sec:2-1:continuumpotato_a}
  \end{subfigure}
  \begin{subfigure}[b]{0.49\textwidth}
    \begin{tikzpicture}[scale=0.8]
    \coordinate (K) at (0,0);
    % Potato
    \draw [fill=black!10,line width=1pt] (K) plot [smooth cycle,tension=0.7] %coordinates {(3,1) (5,1.2) (7,1) (8,3) (7,4.5) (5,4.5) (2,4) (1.7,2.5)};
    coordinates {(3,1) (7,1) (8,3) (7,4.475) (5,4.5) (2,4) (1.7,2.5)};
    \node[ ] at (3.25,2.15) {$\dom$};
    % Crack Surface
    \fill[black, path fading=fade out, draw=none] (5,2.5) circle (0.3);
    \draw[line width=0.1pt,black] (5,0.75) to (5,2.5);
    \draw[line width=0.1pt,black] (5,0.75) to (5,2.5);
    \shade [top color=black,bottom color=black!10,shading angle=90] (5,0.78) rectangle (5.3,2.5);
    \shade [top color=black!10,bottom color=black,shading angle=90] (4.7,0.78) rectangle (5,2.5);
    \draw[->,line width=1pt,black] (4.3,1.7) to (4.85,1.7);
    \draw[<-,line width=1pt,black] (5.15,1.7) to (5.7,1.7);
    \node[ ] at (5.85,1.7) {$\l$};
    %\draw[line width=0.5pt,black] (5,2.5) to (5.5,2.8);
    %\node[ ] at (6,2.9) {$\surfpfreg$};
    % Dirichlet Boundary
    \draw (K) [line width=2.5pt,black] plot [smooth, tension=0.8] coordinates {(3,4.3) (1.7,3.7) (1.7,2.5)};
    \node[ ] at (1.15,4) {$\surfarg{D}{u}$};
    % Neumann Boundary
    %\draw (K) [line width=2.5pt,black] plot [smooth, tension=0.7] coordinates {
    %(7.95,3.4) (7.05,4.5) (5.05,4.5)};
    \node[ ] at (6.75,5.05) {$\surfarg{N}{u}$};
    \end{tikzpicture}
    \caption{diffused (smeared) crack}
    \label{fig:sec:2-1:continuumpotato_b}
  \end{subfigure}
  \caption{A solid, $\dom \in \mathbb{R}^2$, embedded with (a) discrete crack $\surfcrack$ and (b) diffused (smeared) crack, with Dirichlet and Neumann boundaries indicated as $\surfarg{D}{u}$ and $\surfarg{N}{u}$ respectively. Figure reproduced from \cite{Bharali2021}.}
  \label{fig:sec:2:continuumpotato}
\end{figure*}
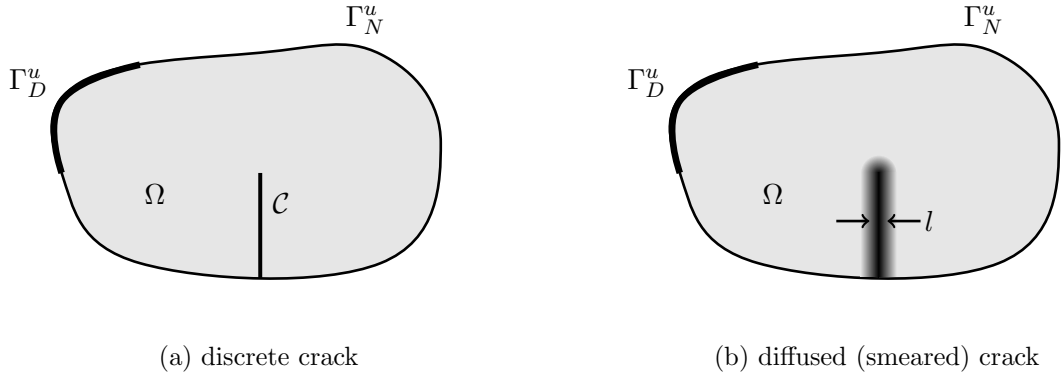

Let $\dom \in \mathbb{R}^{\text{dim}}$ (dim $= 1,2,3$) be a domain occupied by a fracturing continuum. A discrete representation of fracture is shown in Figure \ref{fig:sec:2-1:continuumpotato_a} where the fracture $\mathcal{C}$ may be represented by a cohesive zone fracture interface. Its diffused counterpart, obtained through the phase-field regularisation is presented in Figure \ref{fig:sec:2-1:continuumpotato_b}. Here, the fracture is represented by an auxiliary variable, the phase-field parameter, $\pf \in [0,1]$ within a diffusive (smeared) zone of width $\l>0$, where $\l$ denotes a length scale parameter that controls the width of the diffused zone. The bounds over $\pf$, zero and one indicate the intact material state and total loss of integrity, respectively. Furthermore, the surface $\surf$ of both, the discrete and the diffused fracture continuum is decomposed into a Dirichlet boundary, $\surfarg{D}{u}$ and a Neumann boundary, $\surfarg{N}{u}$, such that $\surf = \surfarg{D}{u} \cup \surfarg{N}{u}$ and $\surfarg{D}{u} \cap \surfarg{N}{u} = \emptyset$.

A general form of the energy functional for the phase-field fracture model, shown in Figure \ref{fig:sec:2-1:continuumpotato_b}, is given by,

\begin{equation}{\label{eqn:sec2:EFunc}}
    \displaystyle E(\bdisp,\pf) = \int_{\dom}^{} g(\pf) \Psi^f(\strain(\bdisp)) \: \text{d}\dom + \int_{\dom}^{} \Psi^r(\strain(\bdisp)) \: \text{d}\dom + \int_{\dom}^{} \dfrac{\gc}{c_w \l} \left( w(\pf) + \l^2 |\bnabla \pf|^2 \right) \: \text{d}\dom,
\end{equation}

\noindent in the absence of any external loading (body and traction forces). Here, $g(\pf)$ is a monotonically decreasing stress-degradation function attached to the fracture driving strain energy $\Psi^f(\strain(\bdisp))$, and $\Psi^r(\strain(\bdisp))$ is the residual strain energy. Moreover, $\gc$ is the Griffith fracture energy, which is a material parameter, and $c_w$ is a normalisation constant associated with the choice of the locally dissipated fracture energy function, $w(\pf)$. Finally, $\strain(\bdisp)$ is the symmetric part of the deformation gradient, where $\bdisp$ is the displacement field. The choice of $w(\pf) = \pf^2$ ($c_w=1/2$) typically denotes a AT2 phase field model, while the choice $w(\pf) = \pf$ ($c_w=8/3$) is often referred to as the AT1 model.

The phase-field model for fracture allows great flexibility in the choice of degradation function $g(\pf)$ and locally dissipated energy function $w(\pf)$, albeit with some restrictions. The degradation function must satify the following criterion: $g(0) = 1$, $g(1) = 0$, and $g'(1) = 0$, to ensure that the fracture driving energy reaches zero for fully developed crack, $\textit{i.e.,}$ for $\pf = 1$. 

Nevertheless, several researchers have proposed different combinations of degradation functions and locally dissipated fracture energy functions, some of which are presented in Table \ref{tab:sec:2:degradationfunctions}. %and \ref{tab:sec:2:localdisspationfunction}.

\begin{table}[ht!]
\centering
\renewcommand{\arraystretch}{1.5} %for vertical spacing of rows
\begin{tabular}{@{}lll@{}}
\toprule
\multicolumn{1}{c}{Type} & \multicolumn{1}{c}{$g(\pf)$} & \multicolumn{1}{c}{Contribution} \\ \midrule
Quadratic & $(1-\pf)^2$ &  Bourdin et. al. \cite{Bourdin2000} \\
Cubic & $s((1-\pf)^3-(1-\pf)^2)+3(1-\pf)^2 - 2(1-\pf)^3$ & Borden et. al. \cite{borden2016} \\ 
Rational & $\dfrac{(1-\pf)^p}{(1-\pf)^p + a_1\pf + a_1 a_2\pf^2 + a_1 a_2 a_3\pf^3}$ & Wu \cite{wu2017} \\ \bottomrule
\end{tabular}
\caption{Stress-degradation functions, popular in the phase-field fracture literature}
\label{tab:sec:2:degradationfunctions}
\end{table}

\begin{comment}
\begin{table}[ht!]
\centering
\renewcommand{\arraystretch}{1.5} %for vertical spacing of rows
\begin{tabular}{@{}lll@{}}
\toprule
\multicolumn{1}{c}{Abbreviated name} & \multicolumn{1}{c}{$w(\pf)$} & \multicolumn{1}{c}{$c_w$} \\ \midrule
AT1 & $\pf$ & $8/3$ \\
AT2 & $\pf^2$ & $2$ \\
PFCZM & $2\pf - \pf^2$ & $\pi$ \\ \bottomrule
\end{tabular}
\caption{Locally dissipated fracture energy functions, popular in the phase-field fracture literature}
\label{tab:sec:2:localdisspationfunction}
\end{table}
\end{comment}

For AT1 and AT2 brittle fracture, the quadratic degradation function proposed in \cite{Bourdin2000} is most commonly adopted. However, it is observed that the AT2 model lacks an initial elastic stage. In order to obtain a linear pre-peak response with the AT2 model, researchers opt for a cubic degradation function proposed in \cite{borden2016}, with $0 < s \leq 1$, determining the slope of $g(\pf)$ in the undamaged state. For quasi-brittle fracture, \cite{wu2017} developed a rational degradation function with model parameters $a_1$, $a_2$, $a_3$, and $p$ is used. With these parameters, the different traction-separation laws are obtained. The reader is referred to \cite{wu2017} for more on this aspect. In this work, the AT2 model is adopted along with quadratic and cubic degradation functions. 

\begin{comment}
Based on their choices, different traction-separation laws can be implemented, as shown in Table \ref{tab:sec:2:TSL_quasi_brittle}.

\begin{table}[ht!]
\centering
\renewcommand{\arraystretch}{1.5} %for vertical spacing of rows
\begin{tabular}{@{}lclll@{}}
\toprule
\multicolumn{1}{c}{Traction-separation law} & \multicolumn{1}{c}{$a_1$} & \multicolumn{1}{c}{$a_2$} & \multicolumn{1}{c}{$a_3$} & \multicolumn{1}{c}{$p$} \\ \midrule
Linear & $\dfrac{4 E_0 \gc}{\pi \l f_t^2}$ & $-0.5$ & $0.0$ & $2.0$ \\
Exponential & " & $2^{5/3}-3$ & $0.0$ & $2.5$ \\
Cornellisen & " & $1.3868$ & $0.9106$ & $2.0$ \\ \bottomrule
\end{tabular}
\caption{Traction-separation laws for quasi-brittle PFCZM}
\label{tab:sec:2:TSL_quasi_brittle}
\end{table}
\end{comment}

Furthermore, the choice of the fracture driving strain energy, $\Psi^f(\strain(\bdisp))$ and the residual, $\Psi^r(\strain(\bdisp))$ is also not unique. Table \ref{tab:sec:2:strainenergysplit} presents a few commonly adopted fracture driving and residual energy. The first model proposed in \cite{Bourdin2000,Bourdin2007} assumes the fracture to driven by the strain energy density $\Psi$, without accounting for tension-compressive asymmetry. Such a model predicts similar fracture in tension and compression. However, most researchers have adopted the notion that fracture is driven by the tensile strain energy density. In this context, \cite{Miehe2010} adopted a spectral decomposition of the strain energy density function. This yielded the tensile strain energy and the compressive strain energies as the fracture driving and residual energies respectively. In an alternative approach, \cite{wu2017} proposed the energy associated with the maximum principal stress $\sigma_1$ as the fracture driving energy. With reference to Table \ref{tab:sec:2:strainenergysplit}, $E$, $\lambda$ and $\mu$ are material constants corresponding to the Young's modulus, Lam\'e constant and shear modulus respectively. The trace operator is given by $tr$, while $\langle \cdot \rangle_{\pm}$ represents the positive/negative Macaulay brackets. Furthermore, $\strain^{\pm}$ indicates the tensile/compressive strain tensors, obtained through spectral decomposition of the strain tensor.

\begin{table}[ht!]
\centering
\renewcommand{\arraystretch}{1.5} %for vertical spacing of rows
\begin{tabular}{@{}llll@{}}
\toprule
\multicolumn{1}{c}{Type} &\multicolumn{1}{c}{$\Psi^f$} & \multicolumn{1}{c}{$\Psi^r$} & \multicolumn{1}{c}{Contribution} \\ \midrule
No Split & $\frac{1}{2} \lambda tr^2(\strain) + \mu \, \strain \colon \strain$ & $0$ & Bourdin et. al. \cite{Bourdin2000,Bourdin2007} \\
Spectral & $\frac{1}{2} \lambda \langle tr(\strain) \rangle_+^2 + \mu \, \strain^+ \colon \strain^+$ & $\frac{1}{2} \lambda \langle tr(\strain) \rangle_-^2 + \mu \, \strain^- \colon \strain^-$ & Miehe et. al. \cite{Miehe2010} \\
Rankine & $\frac{1}{2 E} \langle \sigma_1 \rangle^2_+ $ & - & Wu \cite{wu2017} \\
 \bottomrule
\end{tabular}%
\caption{Strain energy density decompositions in phase-field based fracture analysis.}
\label{tab:sec:2:strainenergysplit}
\end{table}

The additive decomposition of the strain energy densitiy into fracture driving energy and a residual energy renders the displacement sub-problem non-linear. In order to preserve a linear displacement sub-problem, \cite{Ambati2015} proposed a `\textit{hybrid}' approach. With this approach, the degradation function $g(\pf)$ is applied on the entire strain energy density $\Psi$ instead of $\Psi^f$ in the momentum balance equation. As a consequence, the variational consistency of the problem is lost. Nevertheless, the formulation is consistent w.r.t thermodynamic principles. This formulation, referred to as the `hybrid' phase-field fracture model, is adopted in this work.

\subsection{Euler-Lagrange equations}
\label{sec:Euler_Largrange}
The Euler-Lagrange equations for the phase-field model is obtained upon taking the first variation of the energy functional (\ref{eqn:sec2:EFunc}) w.r.t. its solution fields, vector-valued displacements, $\bdisp$ and scalar-valued phase-field, $\pf$. Incorporating the hybrid formulation \cite{Ambati2015}, and with appropriately defined test and trial spaces, the complete problem statement assumes the form:

\smallskip
\begin{problem}\label{Problem1}
Find ($\bdisp$, $\pf$) $\in \mathbb{U} \times \mathbb{P}$ with

\begin{subequations}
\begin{align}
E'(\bdisp,\pf;\testbdisp) & = \int_{\dom}^{} \bigg( g(\pf) \dfrac{\partial \Psi(\strain(\bdisp))}{\partial \strain} \bigg) \colon \strain(\testbdisp)  \: \normalfont\text{d}\dom = 0 \hspace{9.15em}\forall \: \testbdisp \in \mathbb{U}^0, \label{eqn:sec:2:P1MomentumBalance}  \\
E'(\bdisp,\pf;\hat{\pf}) & = \int_{\dom}^{} \bigg( g'(\pf) \Psi^f(\strain(\bdisp)) + \dfrac{\gc}{c_w \, \l} w'(\pf) \bigg) (\hat{\pf} - \pf)  \: \normalfont\text{d}\dom \label{eqn:sec:2:P1PfEvolution} \\ \nonumber
& + \int_{\dom}^{} \dfrac{2 \gc \l}{c_w} \bnabla\pf \cdot \bnabla(\hat{\pf} - \pf)  \: \normalfont\text{d}\dom \geq 0 \hspace{11em} \forall \: \hat{\pf} \in \mathbb{P}. 
\end{align}
\end{subequations}
\noindent considering the Dirichlet boundary conditions $\bdisp^p$ on $\surfarg{D}{u}$ and $\pf^p$ on $\surfarg{D}{\pf}$, and Neumann boundary condition $\trac{u}$ on $\surfarg{N}{u}$. Moreover, the trial and test spaces are given by
\begin{subequations}
\begin{align}
\mathbb{U} & = \{ \bdisp \in [H^1(\dom)]^{\normalfont\text{dim}}| \bdisp = \bdisp^p \text{ on } \surfarg{D}{u} \},  {\label{eqn:sec:2:P1disp_trialspace}} \\ 
\mathbb{U}^0 & = \{ \bdisp \in [H^1(\dom)]^{\normalfont\text{dim}}| \bdisp = \mathbf{0} \text{ on } \surfarg{D}{u} \},  {\label{eqn:sec:2:P1disp_testspace}} \\
\mathbb{P} & = \{ \pf \in [H^1(\dom)] | \pf \geq {}^{n}\pf |\pf = \pf^p \text{ on } \surfarg{D}{\pf} \}.  {\label{eqn:sec2:P1pf_space}}
\end{align}
\end{subequations}
The left superscript $n$ in (\ref{eqn:sec2:P1pf_space}) refers to the previous step in (pseudo) time. {\color{black}\hfill $\blacksquare$}
\end{problem}

Problem \ref{Problem1} belongs to the variational inequality category (see Equation (\ref{eqn:sec:2:P1PfEvolution} and test/trial space (\ref{eqn:sec2:P1pf_space})). The treatment of variational inequality is not new in the phase-field fracture model literature. A review of the different approaches ensuring fracture irreversibility is presented in Section \ref{sec:1}. Adopting the history-variable approach \cite{Miehe2010}, in conjunction with appropriately defined test and trial spaces, the complete problem statement takes the form:

\smallskip
\begin{problem}\label{Problem2}
Find ($\bdisp$, $\pf$) $\in \mathbb{U} \times \mathbb{P}$ with

\begin{subequations}
\begin{align}
E'(\bdisp,\pf;\testbdisp) & = \int_{\dom}^{} \bigg( g(\pf) \dfrac{\partial \Psi(\strain(\bdisp))}{\partial \strain}  \bigg) \colon \strain(\testbdisp)  \: \normalfont\text{d}\dom = 0 \hspace{9.15em}\forall \: \testbdisp \in \mathbb{U}^0, \label{eqn:sec:2-2:P2MomentumBalance} \\
E'(\bdisp,\pf;\testpf) & = \int_{\dom}^{} \bigg( g'(\pf) \mathcal{H} + \dfrac{\gc}{c_w \, \l} w'(\pf) \bigg) \testpf  \: \normalfont\text{d}\dom \label{eqn:sec:2-2:P2PfEvolution} \\ \nonumber
& + \int_{\dom}^{} \dfrac{2 \gc \l}{c_w} \bnabla\pf \cdot \bnabla\testpf  \: \normalfont\text{d}\dom = 0 \hspace{9.15em} \forall \: \testpf \in \mathbb{P}^0,
\end{align}
\end{subequations}
\noindent considering the Dirichlet boundary conditions $\bdisp^p$ on $\surfarg{D}{u}$ and $\pf^p$ on $\surfarg{D}{\pf}$. Moreover, the trial and test spaces are given by
\begin{subequations}
\begin{align}
\mathbb{U} & = \{ \bdisp \in [H^1(\dom)]^{\normalfont\text{dim}}| \bdisp = \bdisp^p \text{ on } \surfarg{D}{u} \},  {\label{eqn:sec:2-2:P2disp_trialspace}} \\ 
\mathbb{U}^0 & = \{ \bdisp \in [H^1(\dom)]^{\normalfont\text{dim}}| \bdisp = \mathbf{0} \text{ on } \surfarg{D}{u} \},  {\label{eqn:sec:2-2:P2disp_testspace}} \\
\mathbb{P} & = \{ \pf \in [H^1(\dom)] | \pf = \pf^p \text{ on } \surfarg{D}{\pf} \},  {\label{eqn:sec2:pf_P2trialspace}} \\
\mathbb{P}^0 & = \{ \pf \in [H^1(\dom)] | \pf = 0 \text{ on } \surfarg{D}{\pf} \}.  {\label{eqn:sec2:P2pf_testspace}}
\end{align}
\end{subequations}
\noindent The history-variable $\mathcal{H}$ is defined as the maximum fracture driving energy $\Psi^f$ over the entire loading history. Mathematically,
\begin{equation}\label{eqn:sec2:P2hist}
\mathcal{H} = \operatorname{max} ({}^{n}\mathcal{H},\Psi^f).
\end{equation}
The left superscript $n$ in (\ref{eqn:sec2:P2hist}) refers to the previous step in (pseudo) time. {\color{black}\hfill $\blacksquare$}
\end{problem}
\smallskip

\begin{remark}
The history-variable approach in Problem \ref{Problem2} results in a variational equality problem, with relaxed test and trial spaces for the phase-field (cf. Problems \ref{Problem1} and \ref{Problem2}).
\end{remark}

%--------------- IGA and Discrete Equations -------------------------------%
\section{Isogeometric analysis and Discrete equations}\label{sec:3}

\subsection{Isogeometric analysis and NURBS}

Isgeometric analysis (IGA) \cite{hughes2005isogeometric} allows an exact representation of complex geometries, such as such as spheres, ellipsoids, paraboloids and hyperboloids. The representation is carried out using polynomial functions, with the Non-Uniform Rational B-splines (NURBS) most commonly adopted. The smoothness of the the NURBS basis functions is advantageous in problems with multi-faceted surface, that can trigger traction oscillations, when simulated using conventional geometry discretization. Furthermore, IGA also offers the ease in obtaining higher-order continuous basis functions with NURBS. As a consequence, it is appealing for higher-order Partial Differential Equations (PDEs). However, NURBS based modeling is often recognized to have significant flaws in constructing watertight geometries using tensor-product meshes. Also, the scale and scope of refining procedures causes the tensor product structure of NURBS to be inefficient leading to erroneous error estimation and improper implementation of adaptivity algorithms. In the context of phase-field fracture model, NURBS-based simulation was carried out in \cite{BORDEN2014100} for the fourth-order phase-field fracture model, albeit without mesh refinement.

The restrictions of the NURBS-based models were mitigated using T-Splines, while keeping the recognizable structure of NURBS algorithms. T-splines alleviate the deficiencies of NURBS by generating a single patch of watertight geometry that can be fine-tuned and coarsened locally. Implementation of T-splines within the framework of IGA has gained a lot of attention. The B\'ezier extraction \cite{scott2011isogeometric} of the basis makes it suitable to be efficient integration into existing finite element programs. However, the linear independence of T-splines is not assured in generic T-meshes. The concept of analysis suitable T-splines was proposed in \cite{buffa2010linear}, which adopts the essential mathematical entities of NURBS, such as linear independence and partition of unity under certain restrictions on the T-mesh, while giving a highly localized and efficient refining algorithm. As an alternative to T-Splines, PHT-splines were proposed, which are generalization of B-Splines over hierarchical T-meshes. The local refinement algorithm for PHT-splines is extremely efficient, and easy to implement.

In this section, the basic concepts of PHT-splines based IGA is discussed, which is then used as a discretization technique to solve the phase-field fracture problem. In one-dimension, the PHT-spline representation takes the form,

\begin{equation}\label{eq:T_mesh}
    \mathbb{T}(\xi) = \sum_{i=1}^{n_{cp}}\mathcal{N}_{i,p}(\xi)\mathbf{P}_i,
\end{equation}

where, $\mathcal{N}_{i,p}(\xi)$ indicates the cubic B-spline basis functions with $C^1$ continuity defined over the knot vector $\Xi$, and $n_{cp}$ is the total number of control points defined over the control mesh used to determine the scaffolding of the geometry. Furthermore, $p$ denotes the order of the polynomial, and $\mathbf{P}_i \in \mathcal{R}^d$ is the set of control points in $d$ dimensions for the B-spline curve with the knot vector $\Xi$. The initial set of knot vectors is denoted by the set of vertices $\Xi^{d}$ corresponding to each spatial direction in the parameter space, $\hat {\Omega} = \left[0,1\right]^d$, and is given by:

\begin{equation}\label{eq:knot_vector}
    \Xi^{d} = \{{\xi_0^d, \xi_0^d},{\xi_1^d, \xi_1^d},{\xi_2^d},\ldots,{\xi_{n_i-1}^d, \xi_{n_i-1}^d},{\xi_{n_i}^d, \xi_{n_i}^d}\}.
\end{equation}

\noindent Here, ${\xi_i^d}\leq{\xi_{i+1}^d}$, and ${\xi_0^d}\leq{\xi_{1}^d} =0$, ${\xi_{n_i-1}^d}\leq{\xi_{n_i}^d}=1$. Moreover, $n_i =n_{cp} +p+1$ represents the number of elements in each parametric direction. The knot vector is uniform when the distance between the consecutive knots are constant. For each interior vertex $\xi_i$, there are two supporting basis functions, $\left[\xi_{i-1},\xi_{i+1}\right]$. For cubic polynomials, one of the distinguishing aspects of PHT-splines is that they maintain $C^1$ continuity, where the start and end knots are repeated $p+1$ times, while the interior knots are repeated only once. By repeatedly inserting vertices, all the knot spans can be obtained at the same refinement level, hence converting a PHT-spline to a B-spline. 

\begin{figure}[ht!]
    \centering
    \includegraphics[width = 0.55\textwidth]{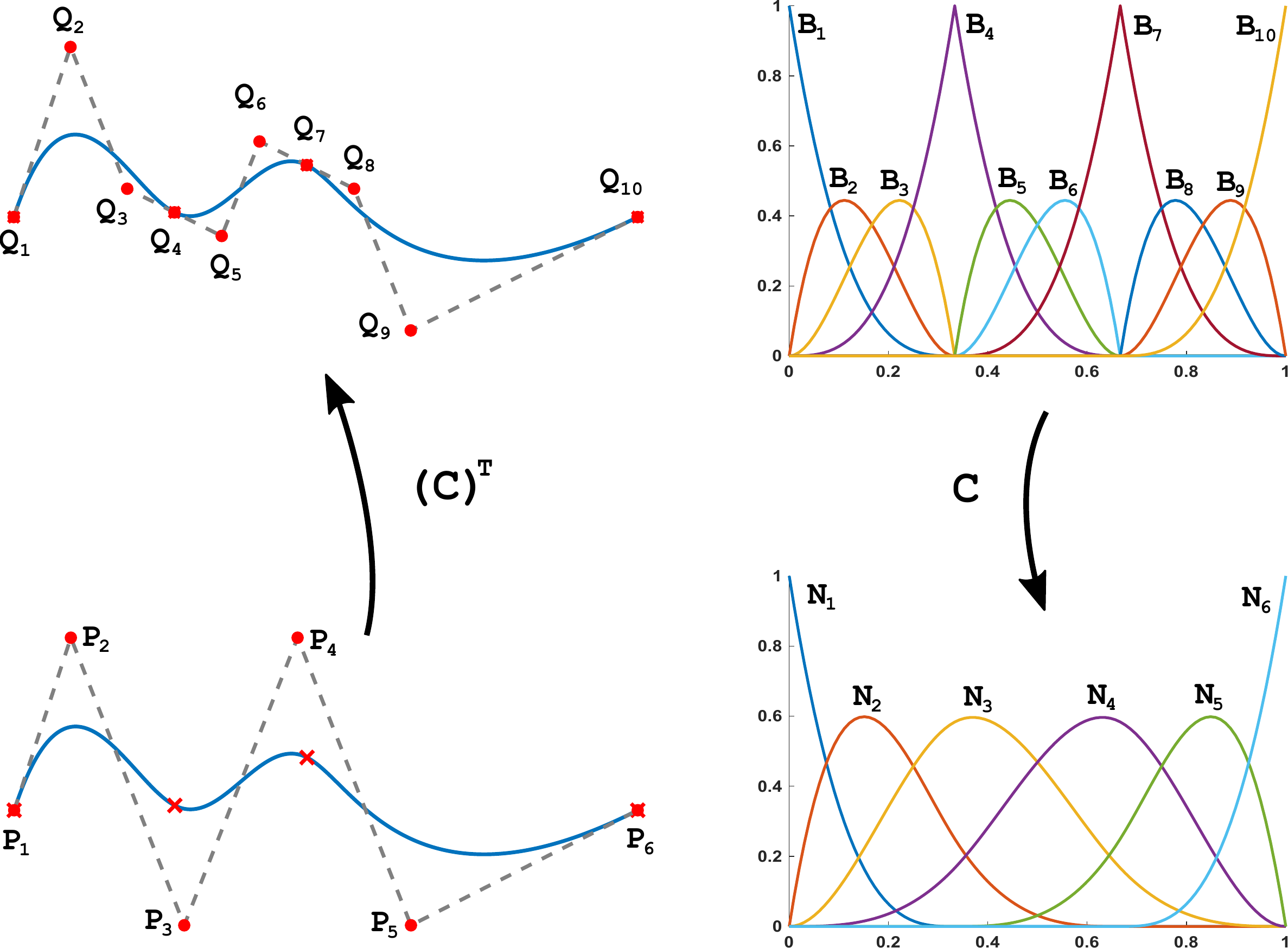}
    \caption{B\'ezier extraction operator implemented for a cubic B-spline curve with $\Xi$ = [0,0,0,0,1/3,2/3,1,1,1,1]. Here, a rational B\'ezier curve, $\mathbf{B}(\bar \xi) = \left\{\beta_{i,p}\right\}_{i=1}^{p+1}$ and its associated control points are $\mathbf{Q} = \left\{Q_i\right\}_{i=1}^{p+1}$ is defined in a reference space $\bar \Omega = [-1,1]$, where $\beta_{i,p}$ are rational Bernstein polynomial basis functions. The B-spline basis functions, $\mathbf{N}(\bar \xi)$ are obtained from the Bernstein basis functions, $\mathbf{B}(\bar \xi)$ using the linear B\'ezier extraction operator $\mathbf{C}$. In case of obtaining the B\'{e}zier control points, the mapping is reversed \cite{goswami2021phase}.}.
  \label{fig:Bez_operator}
\end{figure}

B-splines are non-local, in the sense that, a B-spline typically encompasses more than more element. However, in a finite element framework, a local representation of the B-splines within each element is desired. This local representation of the B-spline is extracted using the B\'ezier decomposition technique. The B\'ezier extraction operator is generated using information from the knot vectors and does not rely on the control points or basis functions. Bernstein polynomials have an edge over NURBS basis functions in terms of implementation because the Bernstein basis functions are the same for all elements, as as observed from Figure \ref{fig:Bez_operator}. Following this idea, the B\'ezier extraction allows for the pre-computation of the Bernstein basis on the reference element. During the simulation, the Bernstein basis function can be mapped to each element, with minimal effort. 

The initial discretization, designated as level 0, is a tensor-product mesh. By splitting the components into $2d$ sub-elements using the cross-insertion approach, the coarse meshes at level $k$ are refined to obtain finer meshes at level $(k+1)$, following the principle of adaptivity. The basis functions for an element on the coarse mesh are replaced by a set of basis functions generated over the refined element in the hierarchical approach. For detailed implementation of the B\'ezier extraction operator and method to obtain the control points on the refined mesh, the reader is referred to \cite{hennig2016bezier}.

\subsection{Discrete phase-field fracture equations}

The IGA proposed for the phase-field fracture model in this work, is similar to the classical Finite Element Analysis (FEA), the only difference being in the basis functions used. For the phase-field fracture model, the Euler-Lagrange equations in Problem \ref{Problem2} is used as the starting point for discretization. Considering the displacement and the phase-field at control points ($\tilde{\bdisp}_i,\tilde{\pf}_i$) as the fundamental unknowns, the corresponding continuous fields ($\bdisp$,$\pf$) are approximated as, 

\begin{equation}\label{eqn:discrete_disp_pf}
  \bdisp = \sum\limits_{i = 1}^m {N_i^{\boldsymbol{u}}\tilde{\boldsymbol{u}}_i}\;\;\;\text{ , }
\pf = \sum\limits_{i = 1}^m N_i^{\varphi}\tilde{\varphi}_i.
\end{equation}

\noindent In the above equation, $N_i^{\boldsymbol{u}}$ and $N_i^{\pf}$ are the basis functions for the displacement and the phase-field, associated with the $i^{\text{th}}$ control point. The total number of control points is given by $m$. The spatial derivatives of the basis functions $N_i^{\boldsymbol{u}}$ and $N_i^{\pf}$ in a two-dimensional case are computed as,

\begin{equation}\label{eqn:Bmatrices}
    \mathbf{B}_i^{\boldsymbol{u}} = \left[ {\begin{array}{*{20}{c}}
  {\begin{array}{*{20}{c}}
  {{N_{i,x}}} \\ 
  0 \\ 
  {{N_{i,y}}}\\
\end{array}}&{\begin{array}{*{20}{c}}
  0 \\ 
  {{N_{i,y}}} \\ 
  {{N_{i,x}}} \\
\end{array}}
\end{array}} \right] \text{ , } 
\mathbf{B}_i^{\pf}  = \left[ {\begin{array}{*{20}{c}}
  {{N_{i,x}}} \\ 
  {{N_{i,y}}}
\end{array}} \right].
\end{equation}

\noindent Here, the subscripts $,x$ and $,y$ indicate spatial derivatives in $x$ and $y$ directions respectively. Using (\ref{eqn:Bmatrices}), the strain $\strain$, and the gradient of the phase-field $\bnabla\pf$ are defined as,

\begin{equation}\label{eqn:strain_gradpf}
  \boldsymbol{\epsilon} = \sum\limits_{i = 1}^m {\mathbf{B}_i^{\boldsymbol{u}} \tilde{\boldsymbol{u}}_i}\;\;\;\text{ , }
\bnabla\pf = \sum\limits_{i = 1}^m {\mathbf{B}_i^{\pf} \tilde{\pf}}.
\end{equation}

\noindent The discrete phase-field fracture problem is obtained upon inserting (\ref{eqn:discrete_disp_pf}-\ref{eqn:strain_gradpf}) in the Euler-Lagrange equations from Problem \ref{Problem2}. Thereafter, (\ref{eqn:sec:2-2:P2MomentumBalance}) and (\ref{eqn:sec:2-2:P2PfEvolution}) are assumed as the internal forces, and stiffness matrix derived from its derivative. This notation is consistent with \cite{de2012nonlinear}. This allows the presentation of the phase-field fracture problem within the incremental iterative framework as:

\smallskip
\begin{dproblem}\label{Problem3}
Compute the solution increment ($\tilde{\Delta\bdisp}$, $\tilde{\Delta\pf}$)$_{i+1}$ in the current iteration $i+1$ using

\begin{subequations}
\begin{equation}\label{dproblem1:discSystem}
\underbrace{
\begin{bmatrix}
\mathbf{K}^{\boldsymbol{uu}} & \mathbf{K}^{\boldsymbol{u}\varphi} \\ 
\mathbf{K}^{\varphi\boldsymbol{u}} & \mathbf{K}^{\varphi\varphi}
\end{bmatrix}_{i}}_{\text{Stiffness matrix}} \begin{Bmatrix}
\Delta \tilde{\bdisp} \\ 
\Delta \tilde{\pf}
\end{Bmatrix}_{i+1}  =  \underbrace{\begin{Bmatrix}
\mathbf{f}^{ext,\boldsymbol{u}} \\ 
\mathbf{f}^{ext,\varphi}
\end{Bmatrix}_{i} - \begin{Bmatrix}
\mathbf{f}^{int,\boldsymbol{u}} \\ 
\mathbf{f}^{int,\varphi}
\end{Bmatrix}_{i}}_{\text{Residual}},
\end{equation}
\noindent and update the solution fields,
\begin{equation}\label{dproblem1:solUpdate}
\begin{Bmatrix}
\tilde{\bdisp} \\ 
\tilde{\pf}
\end{Bmatrix}_{i+1} = \begin{Bmatrix}
\tilde{\bdisp} \\ 
\tilde{\pf}
\end{Bmatrix}_{i} +
\begin{Bmatrix}
\Delta \tilde{\bdisp} \\ 
\Delta \tilde{\pf}
\end{Bmatrix}_{i+1},    
\end{equation}
\noindent until convergence is achieved. The local element stiffness matrices are given by:

\begin{equation}\label{eq:stiff_comp}
\displaystyle
    \begin{split}
        \mathbf{K}^{\boldsymbol{uu}} & = \int_{\Omega} \left[\mathbf{B}^{\boldsymbol{u}}\right]^T \bigg( g(\pf) \dfrac{\partial^2 \Psi}{\partial\strain^2} \bigg) \left[\mathbf{B}^{\boldsymbol{u}}\right]\,d\Omega, \\
        \mathbf{K}^{\boldsymbol{u}\varphi} & =  \int_{\Omega} \left[\mathbf{B}^{\boldsymbol{u}}\right]^T \bigg( g'(\pf) \dfrac{\partial\Psi}{\partial \strain} \bigg) \left[N^{\varphi} \right] \,d\Omega, \\
        \mathbf{K}^{\varphi\boldsymbol{u}} & 
        = \int_{\Omega} \left[N^{\varphi}\right] \bigg( g'(\pf) \dfrac{\partial \mathcal{H}}{\partial \epsilon} \bigg) \left[\mathbf{B}^{\boldsymbol u}\right]\;d\Omega, \\
        \mathbf{K}^{\varphi\varphi} & = 
        \int_{\Omega}\left\{\left[\mathbf{B}^{\varphi}\right]^T \left( \dfrac{2 \gc \l}{c_w} \right) \left[\mathbf{B}^{\varphi}\right] + \left[N^{\varphi}\right] \left(g''(\pf) \mathcal{H} +\dfrac{\gc}{c_w \l} w''(\varphi) \right) \left[N^{\varphi}\right] \right\} \;d\Omega, \\
    \end{split}
\end{equation}
\noindent and the local internal force vectors are computed as
\begin{equation}\label{eq:fint_comp}
    \begin{split}
        \mathbf{f}^{int,\boldsymbol{u}} & = \int_{\Omega} \left[\mathbf{B}^{\boldsymbol{u}}\right]^T \bigg( g(\pf) \dfrac{\partial \Psi}{\partial\strain} \bigg) \,d\Omega,\\
        \mathbf{f}^{int,\varphi} & = 
        \int_{\Omega}\left\{\left[\mathbf{B}^{\varphi}\right]^T \left( \dfrac{2 \gc \l}{c_w} \right) \left[\mathbf{B}^{\varphi}\right] \tilde{\pf} + \left[N^{\varphi}\right]^T \left(g'(\pf) \mathcal{H} +\dfrac{\gc}{c_w \l} w'(\varphi)  \right) \right\} \;d\Omega.
    \end{split}
\end{equation}
The external force vectors $f^{ext,\boldsymbol{u}}$ and $f^{ext,\varphi}$ are considered equal to zero. {\color{black}\hfill $\blacksquare$}
\end{subequations}
\end{dproblem}

%--------------- Monolithic solution technique -------------------------------%
\section{A monolithic solution technique}\label{sec:4}

\subsection{Incremental fracture energy-based arc-length method}

Since the inception of the dissipation-based arc-length solver in \cite{Verhoosel2009dissip}, variants thereof have been utilized for the phase-field based fracture modeling \cite{vignollet2014phase,may2015numerical,singh2016fracture}. Motivated by these studies, in particular, the fracture path controlled path following method proposed in \cite{singh2016fracture}, a fracture energy-based arc-length method is proposed in this work. From the phase-field fracture energy functional (\ref{eqn:sec2:EFunc}), the energy associated with fracture is identified as:
\begin{equation}{\label{eqn:sec2:EFrac}}
    \displaystyle G(\pf) = \int_{\dom}^{} \dfrac{\gc}{c_w \l} \left( w(\pf) + \l^2 |\bnabla \pf|^2 \right) \: \text{d}\dom,
\end{equation}

\noindent and its incremental form is given by

\begin{equation}{\label{eqn:sec2:EFracIncr}}
    \displaystyle \Delta G(\pf,\Delta\pf) = \int_{\dom}^{} \dfrac{\gc}{c_w \l} \left( w'(\pf)\Delta\pf + 2 \l^2 \bnabla \pf \cdot \bnabla \Delta \pf \right) \: \text{d}\dom.
\end{equation}

\noindent Thereafter, an arc-length constraint equation $\Lambda(\pf,\Delta\pf)$ is devised, limiting the incremental phase-field fracture energy in (\ref{eqn:sec2:EFracIncr}) to a certain value $\Delta \tau$. Mathematically,

\begin{equation}\label{eqn:arclength_cons}
    \Lambda^{(\pf,\Delta\pf)} := \displaystyle \Delta G(\pf,\Delta\pf) - \Delta \tau = 0.
\end{equation}

Within a displacement-controlled loading scenario, an additive decomposition of the displacement is carried out as,

\begin{equation}
\tilde{\bdisp} = \mathbf{C} \tilde{\bdisp}_f + \tilde{\bdisp}_p + \tilde{\lambda} \hat{\bdisp}.     
\end{equation}

\noindent Here, $\mathbf{C}$ is a constraint matrix \cite{Shephard1984linear}, $\tilde{\bdisp}_f$, $\tilde{\bdisp}_p$ and $\hat{\bdisp}_p$ are free, prescribed and unit displacements. The load level $\tilde{\lambda}$ acts only on the prescribed displacements $\hat{\bdisp}_p$. With this setup, the discrete problem assumes the form:

\begin{dproblem}\label{Problem4}
Compute the solution increment ($\tilde{\Delta\bdisp}$, $\tilde{\Delta\pf}$, $\Delta\hat{\lambda}$)$_{i+1}$ in the current iteration $i+1$ using

\begin{subequations}
\begin{equation}
\underbrace{
\begin{bmatrix}
\mathbf{K}^{\boldsymbol{uu}} & \mathbf{K}^{\boldsymbol{u}\varphi} & \mathbf{K}^{\boldsymbol{u}\lambda}  \\ 
\mathbf{K}^{\varphi\boldsymbol{u}} & \mathbf{K}^{\varphi\varphi} & \mathbf{K}^{\varphi\lambda}  \\
\mathbf{0} & \mathbf{K}^{\lambda\varphi} & \mathbf{0}
\end{bmatrix}_{i}}_{\text{Stiffness matrix}} \begin{Bmatrix}
\Delta \tilde{\bdisp} \\ 
\Delta \tilde{\pf} \\
\Delta \tilde{\lambda}
\end{Bmatrix}_{i+1}  = - \underbrace{\begin{Bmatrix}
\mathbf{f}^{int,\boldsymbol{u}} \\ 
\mathbf{f}^{int,\varphi} \\
\Lambda^{(\varphi,\Delta\varphi)}
\end{Bmatrix}_{i}}_{\text{Residual}}
\end{equation}
\noindent and update the solution fields,
\begin{equation}\label{eq:arcL_solUpdate}
\begin{Bmatrix}
\tilde{\bdisp} \\ 
\tilde{\pf} \\
\tilde{\lambda}
\end{Bmatrix}_{i+1} = \begin{Bmatrix}
\tilde{\bdisp} \\ 
\tilde{\pf}  \\
\tilde{\lambda}
\end{Bmatrix}_{i} +
\begin{Bmatrix}
\Delta \tilde{\bdisp} \\ 
\Delta \tilde{\pf}  \\
\Delta \tilde{\lambda}
\end{Bmatrix}_{i+1}.    
\end{equation}
\noindent until convergence is achieved. The local element stiffness matrices $\mathbf{K}^{\boldsymbol{uu}}$, $\mathbf{K}^{\boldsymbol{u}\varphi}$, $\mathbf{K}^{\varphi\boldsymbol{u}}$ and $\mathbf{K}^{\varphi\varphi}$, and the local internal force vectors $\mathbf{f}^{int,\boldsymbol{u}}$ and $\mathbf{f}^{int,\varphi}$ remain same as that presented in Problem \ref{Problem3}. The additional matrices are given by

\begin{equation}\label{eq:addl_stiff_comp}
\displaystyle
    \begin{split}
        \mathbf{K}^{\boldsymbol{u}\lambda} & = \mathbf{K}^{\boldsymbol{uu}} \: \hat{\bdisp}, \\
        \mathbf{K}^{\varphi\lambda} & = \mathbf{K}^{\varphi\boldsymbol{u}} \: \hat{\bdisp}, \\
        \mathbf{K}^{\lambda\varphi} & = \dfrac{\gc}{c_w \l}
        \int_{\Omega}\left\{ \left[N^{\varphi}\right]^T \big( w''(\pf)\Delta\pf + w'(\pf) \big) + 2 \l^2 \left[\mathbf{B}^{\varphi}\right]^T \cdot \big( \bnabla\Delta\pf + \bnabla\pf  \big)  \right\} \;d\Omega,
    \end{split}
\end{equation}
\noindent and the additional internal force vector is computed as
\begin{equation}\label{eq:addl_fint_comp}
    \begin{split}
        \Lambda^{(\varphi,\Delta\varphi)} & = \int_{\dom}^{} \dfrac{\gc}{c_w \l} \left( w'(\pf)\Delta\pf + 2 \l^2 \bnabla \pf \cdot \bnabla \Delta \pf \right) \: \text{d}\dom - \Delta\tau.
    \end{split}
\end{equation}
\end{subequations}
{\color{black}\hfill $\blacksquare$}
\end{dproblem}

In this work, the phase-field fracture problem is solved using both the conventional displacement-controlled solution scheme (Problem \ref{Problem3}) and the Arc-length method (Problem \ref{Problem4}). Algorithm \ref{alg:dispToArc} presented in Appendix \ref{appA:algorithms} combines both within a single monolithic solution strategy. The time-stepping commences with a displacement-controlled solution scheme and the Newton-Raphson method. Upon convergence, the incremental dissipation $\Delta G$ is computed using (\ref{eqn:sec2:EFracIncr}). If $\Delta G$ is greater than a certain user-defined switch energy, the method is switched to the Arc-length method. Additionally, $\Delta \tau$ is set to the switch energy, and $\Delta \lambda$ is set to zero. Upon achieving convergence in a certain step, an increment in $\Delta \lambda$ is carried out, subjected to a maximum value of $\Delta\tau_{max}$.

\subsection{Adaptive under-relaxation scheme}

An under-relaxation scheme introduces a scalar parameter $\beta \in (0,1)$ such that the solution field is updated using,

\begin{equation}\label{eq:underRelaxSolUpdate}
\begin{Bmatrix}
\tilde{\bdisp} \\ 
\tilde{\pf} \\
\tilde{\lambda}
\end{Bmatrix}_{i+1} = \begin{Bmatrix}
\tilde{\bdisp} \\ 
\tilde{\pf}  \\
\tilde{\lambda}
\end{Bmatrix}_{i} + \beta
\begin{Bmatrix}
\Delta \tilde{\bdisp} \\ 
\Delta \tilde{\pf}  \\
\Delta \tilde{\lambda}
\end{Bmatrix}_{i+1}.        
\end{equation}

\noindent When $\beta$ is set to one, the Newton-Raphson method is recovered, otherwise, the method maybe referred to as a modified Newton approach. Under-relaxation schemes are usually robust, however, it may reduce the rate of convergence of the problem \cite{STORVIK2021113822}. For the phase-field fracture problem, the under-relaxation scheme is adopted to prevent divergence due to an ill behaving stiffness matrix.

Algorithm \ref{alg:dispToArcUnderRelax} in Appendix \ref{appA:algorithms} presents the under-relaxation adopted in this work. In this scheme, $\beta$ starts with a value one, corresponding to a full Newton-Raphson update. When the convergence is not achieved, the value to $\beta$ is reduced by a factor 1.25. This reduction is carried out twice before performing a reduction in the prescribed incremental dissipation $\Delta\tau$. The motivation behind this is to try the current dissipation step with smaller solution increments within the iterative process of the Newton-Raphson method. Such an approach prevents divergence due to an ill behaving stiffness matrix.

\subsection{Adaptive mesh refinement and solution transfer}

The phase-field fracture model requires an extremely fine mesh to resolve the crack zone in the computational domain, $\dom$. A sharp crack topology is recovered in the limit as $\l \to 0$ \cite{Bourdin2000}. With fracture length-scales $\l$ very small compared to $|\dom|$, an uniformly refined mesh enormously increases the resources required in terms of computing power and storage. In this work, the novel monolithic solver is integrated with an efficient adaptive mesh refinement (AMR) scheme. The elements of the mesh are chosen for refinement based on a critical threshold value of $\pf$, $\pf_{\text{threshold}}$ \cite{Heister2015,goswami2019adaptive}, which is typically referred as the the refinement indicator. To locally refine the crack path, polynomial splines over hierarchical T-meshes (PHT-splines) are used within the framework of IGA. The PHT-splines possess a very efficient and easy to implement local refinement algorithm. The hierarchical approach replaces the basis functions for an element on the coarse mesh with a set of basis functions constructed over the refined element. The refinement of an element originally defined on the coarse mesh is restricted by a pre-decided maximum number of allowable refinements, to avoid repeated refinements of the elements which are already in the cracked domain.

An adaptive h-refinement scheme is adopted in this work, in which the order of the basis functions remains constant throughout the refinement process. Within the simulation, a series of hierarchical meshes evolve. The mesh during the onset of the simulation is denoted as the base mesh or Level 0 mesh. At any hierarchical level, say `k', some (parent) elements are marked for refinement, following which they are sub-divided into $2^{\text{dim}}$ (children) elements. Once a parent element is refined, it becomes inactive and its children take the place in the computational domain as active elements. Finally, for computational efficiency, the basis functions are computed only for the children elements upon refinement, and not for the entire computational domain.

The field variables are projected from the coarser mesh to the finer mesh for each refined element, repeatedly during the mesh refinement. To avoid re-computation of the problem from the beginning, a variable transfer is required. The discretized variables include the field variables such as $\bdisp$ and $\pf$, computed at the control points that are required to be transferred to the new element. The projection of the field variables from a coarser mesh to a finer mesh is implemented using a similar technique to that described in \cite{hennig2016bezier}. For computing the new control points, instead of projecting the geometrical information at the basis vertex, we project $\pf$ on the finer mesh.

%\subsection{Choice of maximum dissipation energy, $\Delta\tau_{max}$}

%The maximum dissipation energy $\Delta\tau_{max}$, allowed in a single step plays a key role in the fracture energy-based arc-length method, and more so with the use of an adaptive mesh refinement technique. The energy dissipated in a single step should not be too large, so that the crack propagates across several elements. The reason for this is that the solution algorithm may fail when the elements are not sufficiently refined. In this work, a simple initial estimate of $\Delta\tau_{max}$ is proposed, equating it to the locally dissipated fracture energy function in (\ref{eqn:sec2:EFunc}) integrated over an element ($e$). Mathematically, 

%\begin{equation}{\label{eqn:sec4:dTauMax}}
%    \displaystyle \Delta\tau_{max} = \int_{\dom,e}^{} \dfrac{\gc}{c_w \l} w(\pf)  \: \text{d}\dom,
%\end{equation}

%--------------- Numerical experiments -------------------------------%
\section{Numerical experiments}\label{sec:5}

In this section, numerical experiments are carried out on representative fracture problems. These include a tapered bar under tension, a specimen with an eccentric hole under tension, a single edge notched specimen under tensile and later shear loading. For all problems, the geometry, material properties as well as loading conditions are presented in the respective sub-sections. The load-displacement curves and the phase-field fracture topology at the final step of the analysis are also presented therein.

A residual-based convergence criterion is adopted in this work. More specifically, the iterations pertaining to the Newton-Raphson method is terminated, when the norm of the residual is less than $1e-3$.

%\SG{May be a nice idea would be to mention critical theshold value for each examples.}

\subsection{Tapered Bar under Tension (TBT)}\label{sec5:TBtension}

The first example in the numerical study section is a tapered bar under tensile loading, as shown in Figure \ref{sec5:fig:TPTdiagram}. The bar has dimensions $5$ [mm] in length, $0.75$ [mm] and $2$ [mm] width of fixed end and the prescribed loading edges respectively. The loading in applied in the form of prescribed displacement increment $\tilde{\lambda}\hat{\bdisp}$, where $\hat{\bdisp}$ is a unit load vector and $\tilde{\lambda}$ is the load factor. When the analysis is started, the displacement-controlled approach is adopted and $\hat{\lambda}$ is incremented in steps of $1e-2$ [mm]. Following the switch to the arc-length method, $\hat{\lambda}$ becomes an unknown variable, and is solved with the arc-length constraint equation. The additional parameters required for the analysis are presented in Table \ref{sec5:table:TBTparams}.

\begin{figure}[ht]
\begin{minipage}[b]{0.45\linewidth}
\centering
\begin{tikzpicture}[scale=0.5]
    \coordinate (K) at (0,0);
    % Taper
    \draw[line width=0.75pt,black] (0,0) to (0,-0.75);
    \draw[line width=0.75pt,black] (0,-0.75) to (10,-2);
    \draw[line width=0.75pt,black] (10,-2) to (10,2);
    \draw[line width=0.75pt,black] (10,2) to (0,0.75);
    \draw[line width=0.75pt,black] (0,0.75) to (0,0);
    % Right Boundary
    %\draw[line width=0.75pt,black] (2.5,2.75) to (-2.5,2.75);
    \draw[->,line width=1.5pt,black] (10.2,0.0) to (10.8,0.0);
    \node[ ] at (11.5,0.0) {$\tilde{\lambda}\hat{\bdisp}$};
    % Left edge
    \draw[line width=0.75pt,black] (0,-1.25) to (0,1.25);
    % Left boundary
    \draw[line width=0.75pt,black] (0,1.25) to (-0.25,1.0);
    \draw[line width=0.75pt,black] (0,1.00) to (-0.25,0.75);
    \draw[line width=0.75pt,black] (0,0.75) to (-0.25,0.50);
    \draw[line width=0.75pt,black] (0,0.50) to (-0.25,0.25);
    \draw[line width=0.75pt,black] (0,0.25) to (-0.25,0.00);
    \draw[line width=0.75pt,black] (0,0.00) to (-0.25,-0.25);
    \draw[line width=0.75pt,black] (0,-0.25) to (-0.25,-0.50);
    \draw[line width=0.75pt,black] (0,-0.50) to (-0.25,-0.75);
    \draw[line width=0.75pt,black] (0,-0.75) to (-0.25,-1.00);
    \draw[line width=0.75pt,black] (0,-1.00) to (-0.25,-1.25);
    \draw[line width=0.75pt,black] (0,-1.25) to (-0.25,-1.50);
    \end{tikzpicture}
\caption{TBT experiment}
\label{sec5:fig:TPTdiagram}
\end{minipage}
\begin{minipage}[b]{0.45\linewidth}
\centering
\begin{tabular}{ll} \hline
  \textbf{Parameters} & \textbf{Value} \\ \hline
  Model & AT2 \\
  $\Psi^f$ & No Split \\
  $\lambda$ & 0.0 [MPa] \\
  $\mu$ & 50.0 [MPa] \\
  $\gc$ & 1.0 [N/mm] \\
  $\l$ & 0.25 [mm] \\
  $\pf_{\text{threshold}}$ & 0.2 \\  
  $\Delta\tau_{max}$ & 0.0125 [N] \\ \hline
  \end{tabular}
\captionof{table}{Model parameters}
\label{sec5:table:TBTparams}
\end{minipage}
\end{figure}

\begin{figure}[!ht]
  \begin{subfigure}[t]{0.45\textwidth}
  \centering
    \begin{tikzpicture}[thick,scale=0.95, every node/.style={scale=0.95}]
    \begin{axis}[legend style={draw=none}, legend columns = 2,
      transpose legend, ylabel={Load\:[N]},xlabel={Displacement\:[mm]}, xmin=0, ymin=0, xmax=0.45, ymax=25, yticklabel style={/pgf/number format/.cd,fixed,precision=2},
                 every axis plot/.append style={very thick}]
    \pgfplotstableread[col sep = comma]{./Data/TaperBar/lodi_quad.txt}\Adata;
    \pgfplotstableread[col sep = comma]{./Data/TaperBar/lodi_cubics0_01.txt}\Bdata;
    \pgfplotstableread[col sep = comma]{./Data/TaperBar/lodi_cubics0_1.txt}\Cdata;
    \pgfplotstableread[col sep = comma]{./Data/TaperBar/lodi_cubics1.txt}\Ddata;
    \addplot [ 
           color=black, 
%           only marks, 
           mark=*, 
           mark size=0.15pt, 
         ]
         table
         [
           x expr=\thisrowno{1}, 
           y expr=\thisrowno{0}
         ] {\Adata};
         \addlegendentry{Quadratic}
%    \addplot [ 
%           color=red, 
%           only marks, 
%           mark=*, 
%           mark size=0.75pt, 
%         ]
%         table
%         [
%           x expr=\thisrowno{1}, 
%           y expr=\thisrowno{0}
%         ] {\Bdata};
%         \addlegendentry{Cubic: $s = 0.01$}
    \addplot [ 
           color=red, 
%           only marks, 
           mark=*, 
           mark size=0.5pt, 
         ]
         table
         [
           x expr=\thisrowno{1}, 
           y expr=\thisrowno{0}
         ] {\Cdata};
         \addlegendentry{Cubic: $s = 0.1$}   
    \addplot [ 
           color=blue, 
%           only marks, 
           mark=*, 
           mark size=0.25pt, 
         ]
         table
         [
           x expr=\thisrowno{1}, 
           y expr=\thisrowno{0}
         ] {\Ddata};
         \addlegendentry{Cubic: $s = 1.0$}      
    \end{axis}
    \end{tikzpicture}
    \caption{ }
    \label{sec5:fig:TBT_lodi}
  \end{subfigure}
  \hfill
  \begin{subfigure}[t]{0.45\textwidth}
  \centering
    \begin{tikzpicture}
    \node[inner sep=0pt] () at (0,0)
    {\includegraphics[width=5.5cm,trim=8.6cm 2.38cm 2.65cm 7.15cm, clip]{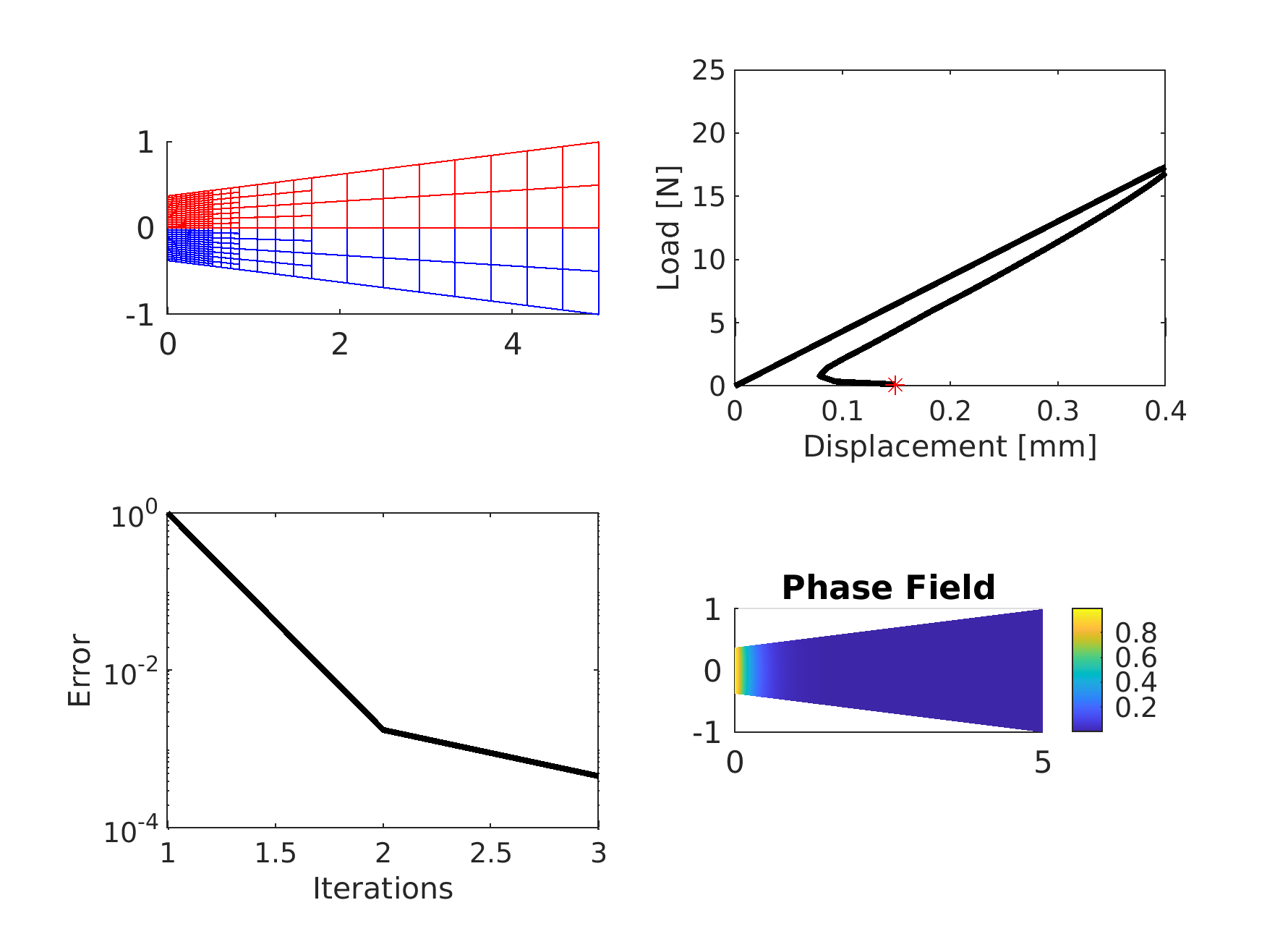}};
    \node[inner sep=0pt] () at (-1.3,-3.05)
    {\begin{axis}[
    hide axis,
    scale only axis,
    height=0pt,
    width=0pt,
    %colormap/bluered,
    colorbar horizontal,
    point meta min=0,
    point meta max=1,
    colorbar style={
        width=4.85cm,
        xtick={0,0.5,1.0},
        xticklabel style = {yshift=-0.075cm}
    }]
    \addplot [draw=none] coordinates {(0,0)};
    \end{axis}};
    \node[inner sep=0pt] () at (0,-3.75) {$\pf$};
    \end{tikzpicture}
    \caption{ }
    \label{sec5:fig:TBT_failure}
  \end{subfigure}
  \caption{Figure (a) presents the load-displacement curves for the tapered bar under tension. The legend entries correspond to the choice of degradation functions. Figure (b) shows the distribution of the phase-field variable at the final step of the analysis.}
\end{figure}

Figure \ref{sec5:fig:TBT_lodi} presents the load-displacement curves for the tapered bar under tension. For both, quadratic and cubic degradation function, the specimen exhibits a snap-back behaviour upon localisation. Furthermore, consistent with other studies in the literature, the cubic degradation function demonstrates a linear pre-peak behaviour for small values of $s$ ($<1$). It is also observed that the use of the cubic degradation function yields a higher peak load compared to the quadratic degradation function. The reason for this behaviour can be explained with analytical studies on 1D bar \cite{pham2011}. Finally, in Figure \ref{sec5:fig:TBT_failure}, the phase-field fracture topology at the final step of the analysis is presented, where the fracture is seen to appear on the fixed end.

\subsection{Specimen with an eccentric hole under tension (EH)}

A unit square (in mm) embedded with a eccentric hole \cite{May2016}, as shown in Figure \ref{sec5:fig:EHdiagram}, is considered for numerical study. The hole has a radius 0.2 [mm] and is centered at (0.6,0.0). A quasi-static loading is applied at the top boundary in the form of prescribed displacement increment $\tilde{\lambda}\hat{\bdisp}$, where $\hat{\bdisp}$ is a unit load vector and $\tilde{\lambda}$ is the load factor. The loading in applied in the form of prescribed displacement increment $\tilde{\lambda}\hat{\bdisp}$, where $\hat{\bdisp}$ is a unit load vector and $\tilde{\lambda}$ is the load factor. When the analysis is started, the displacement-control approach is adopted and $\hat{\lambda}$ is incremented in steps of $1e-4$ [mm]. Following the switch to the arc-length method, $\hat{\lambda}$ becomes an unknown variable, and is solved with the arc-length constraint equation. The bottom boundary remains fixed. The additional parameters required for the analysis are presented in Table \ref{sec5:table:EHparams}.

\begin{figure}[ht]
\begin{minipage}[b]{0.45\linewidth}
\centering
\begin{tikzpicture}[scale=0.7]
    \coordinate (K) at (0,0);
    % Square
    \draw [fill=black!20] (-2.5,-2.5) rectangle (2.5,2.5);
    \draw [fill=white] (0.5,0) circle (1.0);
    %\draw[line width=0.75pt,black] (-2.5,-2.5) to (2.5,-2.5);
    %\draw[line width=0.75pt,black] (2.5,-2.5) to (2.5,2.5);
    %\draw[line width=0.75pt,black] (2.5,2.5) to (-2.5,2.5);
    %\draw[line width=0.75pt,black] (-2.5,2.5) to (-2.5,-2.5);
    %\draw[line width=1.5pt,red] (-2.5,0) to (0,0);
    % Top Boundary
    \draw[line width=0.75pt,black] (2.5,2.75) to (-2.5,2.75);
    \draw[->,line width=1.5pt,black] (0.0,3.0) to (0.0,3.5);
    \node[ ] at (0,3.85) {$\tilde{\lambda}\hat{\bdisp}$};
    % Bottom Boundary
    \draw[line width=0.75pt,black] (-2.5,-2.5) to (-2.75,-2.75);
    \draw[line width=0.75pt,black] (-2.25,-2.5) to (-2.5,-2.75);
    \draw[line width=0.75pt,black] (-2.0,-2.5) to (-2.25,-2.75);
    \draw[line width=0.75pt,black] (-1.75,-2.5) to (-2.0,-2.75);
    \draw[line width=0.75pt,black] (-1.5,-2.5) to (-1.75,-2.75);
    \draw[line width=0.75pt,black] (-1.25,-2.5) to (-1.5,-2.75);
    \draw[line width=0.75pt,black] (-1.0,-2.5) to (-1.25,-2.75);
    \draw[line width=0.75pt,black] (-0.75,-2.5) to (-1.0,-2.75);
    \draw[line width=0.75pt,black] (-0.5,-2.5) to (-0.75,-2.75);
    \draw[line width=0.75pt,black] (-0.25,-2.5) to (-0.5,-2.75);
    \draw[line width=0.75pt,black] (0.0,-2.5) to (-0.25,-2.75);
    \draw[line width=0.75pt,black] (0.25,-2.5) to (0.0,-2.75);
    \draw[line width=0.75pt,black] (0.5,-2.5) to (0.25,-2.75);
    \draw[line width=0.75pt,black] (0.75,-2.5) to (0.5,-2.75);
    \draw[line width=0.75pt,black] (1.0,-2.5) to (0.75,-2.75);
    \draw[line width=0.75pt,black] (1.25,-2.5) to (1.0,-2.75);
    \draw[line width=0.75pt,black] (1.5,-2.5) to (1.25,-2.75);
    \draw[line width=0.75pt,black] (1.75,-2.5) to (1.5,-2.75);
    \draw[line width=0.75pt,black] (2.0,-2.5) to (1.75,-2.75);
    \draw[line width=0.75pt,black] (2.25,-2.5) to (2.0,-2.75);
    \draw[line width=0.75pt,black] (2.5,-2.5) to (2.25,-2.75);
    % \node[ ] at (-0.5,-2.05) {$\tilde{\bdisp}_{p} = 0$};
    \end{tikzpicture}
\caption{EH experiment}
\label{sec5:fig:EHdiagram}
\end{minipage}
\begin{minipage}[b]{0.45\linewidth}
\centering
\begin{tabular}{ll} \hline
  \textbf{Parameters} & \textbf{Value} \\ \hline 
  Model & AT2 \\
  $\Psi^f$ & Rankine \\
  $\lambda$ & 121.154 [GPa] \\
  $\mu$ & 80.769 [GPa] \\
  $\gc$ & 2700 [N/m] \\
  $\l$ & 2e-2 [mm] \\
  $\pf_{\text{threshold}}$ & 0.2 \\  
  $\Delta\tau_{max}$ & 0.05 [N/mm$^2$] \\ \hline
  \end{tabular}
\captionof{table}{Model parameters}
\label{sec5:table:EHparams}
\end{minipage}
\end{figure}

\begin{figure}[!ht]
 \begin{subfigure}[t]{0.45\textwidth}
  \centering
    \begin{tikzpicture}[thick,scale=0.95, every node/.style={scale=0.95}]
    \begin{axis}[legend style={draw=none}, legend columns = 2,
      transpose legend, ylabel={Load\:[N]},xlabel={Displacement\:[mm]}, xmin=0, ymin=0, xmax=0.008, ymax=1600, yticklabel style={/pgf/number format/.cd,fixed,precision=2},
                 every axis plot/.append style={very thick}]
    \pgfplotstableread[col sep = comma]{./Data/EH/lodi_quad.txt}\Adata;
    %\pgfplotstableread[col sep = comma]{./Data/SENT/lodi_cubics0_01.txt}\Bdata;
    \pgfplotstableread[col sep = comma]{./Data/EH/lodi_cubics0_1.txt}\Cdata;
    \pgfplotstableread[col sep = comma]{./Data/EH/lodi_cubics1.txt}\Ddata;
    \addplot [ 
           color=black, 
%           only marks, 
           mark=*, 
           mark size=0.25pt, 
         ]
         table
         [
           x expr=\thisrowno{1}, 
           y expr=\thisrowno{0}
         ] {\Adata};
         \addlegendentry{Quadratic}
%    \addplot [ 
%           color=red, 
%           only marks, 
%           mark=*, 
%           mark size=0.25pt, 
%         ]
%         table
%         [
%           x expr=\thisrowno{1}, 
%           y expr=\thisrowno{0}
%         ] {\Bdata};
%         \addlegendentry{Cubic: $s = 0.01$}     
    \addplot [ 
           color=red, 
%           only marks, 
           mark=*, 
           mark size=0.25pt, 
         ]
         table
         [
           x expr=\thisrowno{1}, 
           y expr=\thisrowno{0}
         ] {\Cdata};
         \addlegendentry{Cubic: $s = 0.1$}     
    \addplot [ 
           color=blue, 
%           only marks, 
           mark=*, 
           mark size=0.25pt, 
         ]
         table
         [
           x expr=\thisrowno{1}, 
           y expr=\thisrowno{0}
         ] {\Ddata};
         \addlegendentry{Cubic: $s = 1.0$}          
    \end{axis}
    \end{tikzpicture}
    \caption{ }
    \label{sec5:fig:EH_lodi}
  \end{subfigure}
  \hfill
  \begin{subfigure}[t]{0.45\textwidth}
  \centering
    \begin{tikzpicture}
    \node[inner sep=0pt] () at (0,0)
    {\includegraphics[width=5.5cm,trim=8.35cm 1.5cm 2.65cm 6.25cm, clip]{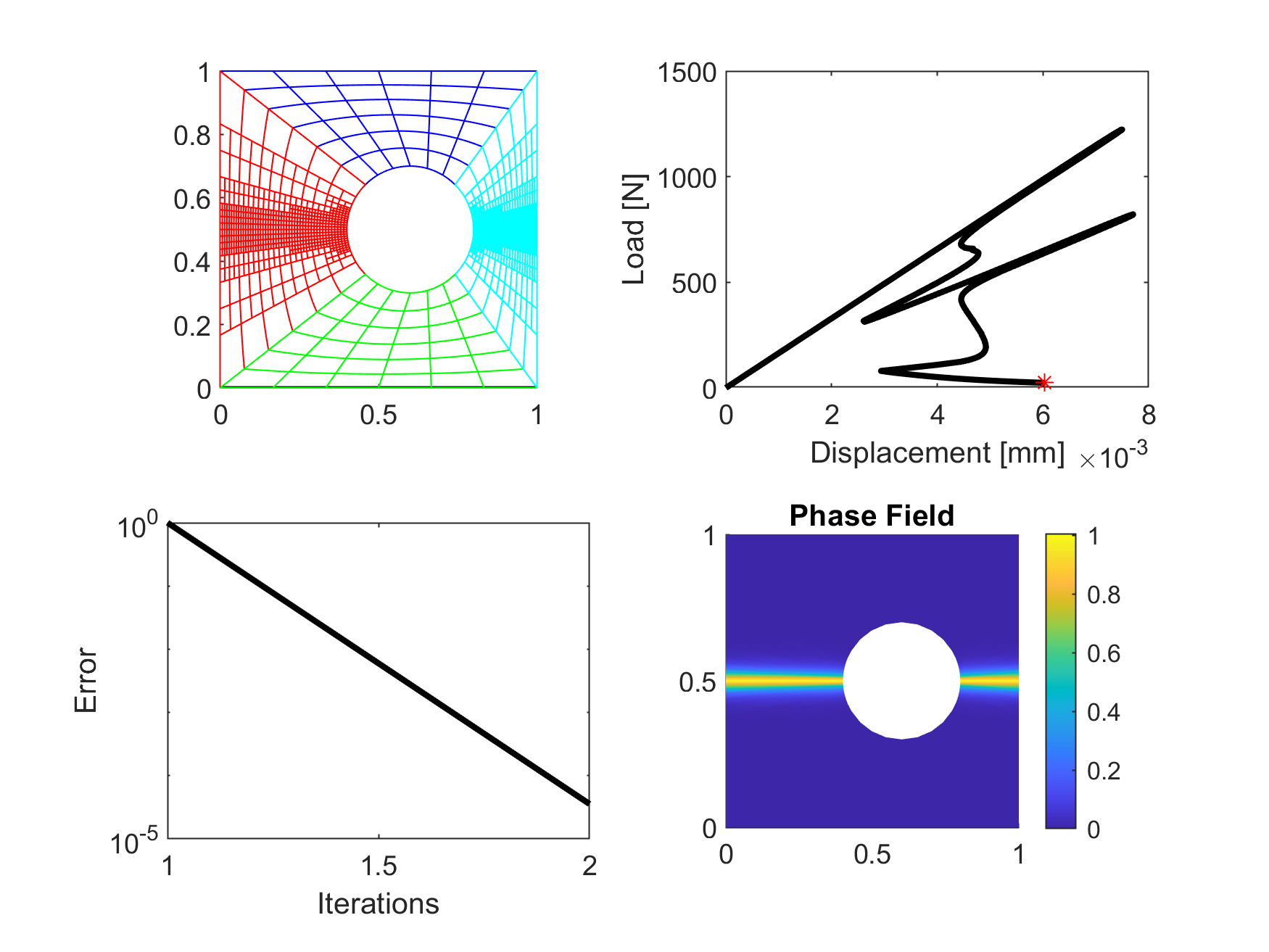}};
    \node[inner sep=0pt] () at (-1.3,-3.05)
    {\begin{axis}[
    hide axis,
    scale only axis,
    height=0pt,
    width=0pt,
    %colormap/bluered,
    colorbar horizontal,
    point meta min=0,
    point meta max=1,
    colorbar style={
        width=4.85cm,
        xtick={0,0.5,1.0},
        xticklabel style = {yshift=-0.075cm}
    }]
    \addplot [draw=none] coordinates {(0,0)};
    \end{axis}};
    \node[inner sep=0pt] () at (0,-3.75) {$\pf$};
    \end{tikzpicture}
    \caption{ }
    \label{sec5:fig:EH_failure}
  \end{subfigure}
  \caption{Figure (a) presents the load-displacement curves for the single edge notched specimen under shear. The legend entries correspond to the choice of degradation functions. Figure (b) shows the distribution of the phase-field variable at the final step of the analysis.}
\end{figure}

\begin{figure}[!ht]
 \begin{subfigure}[t]{0.22\textwidth}
  \centering
    \begin{tikzpicture}[scale=0.5]
    \node[inner sep=0pt] () at (0,0)
    {\includegraphics[width=3.0cm,trim=8.35cm 1.5cm 2.65cm 6.25cm, clip]{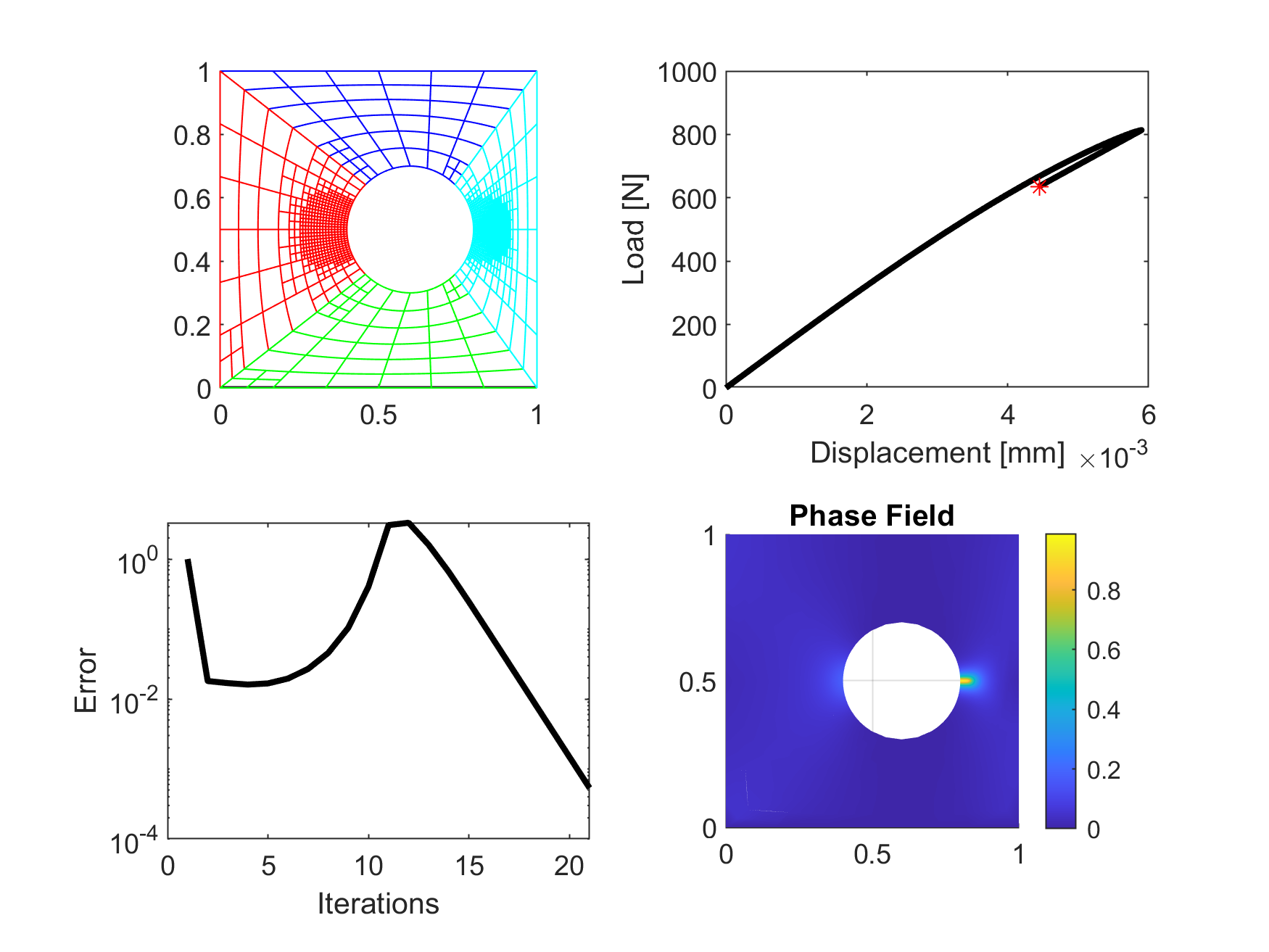}};
    \end{tikzpicture}
    \caption{ }
  \end{subfigure}
  %\hfill
  %
  \begin{subfigure}[t]{0.22\textwidth}
  \centering
    \begin{tikzpicture}[scale=0.5]
    \node[inner sep=0pt] () at (0,0)
    {\includegraphics[width=3.0cm,trim=8.35cm 1.5cm 2.65cm 6.25cm, clip]{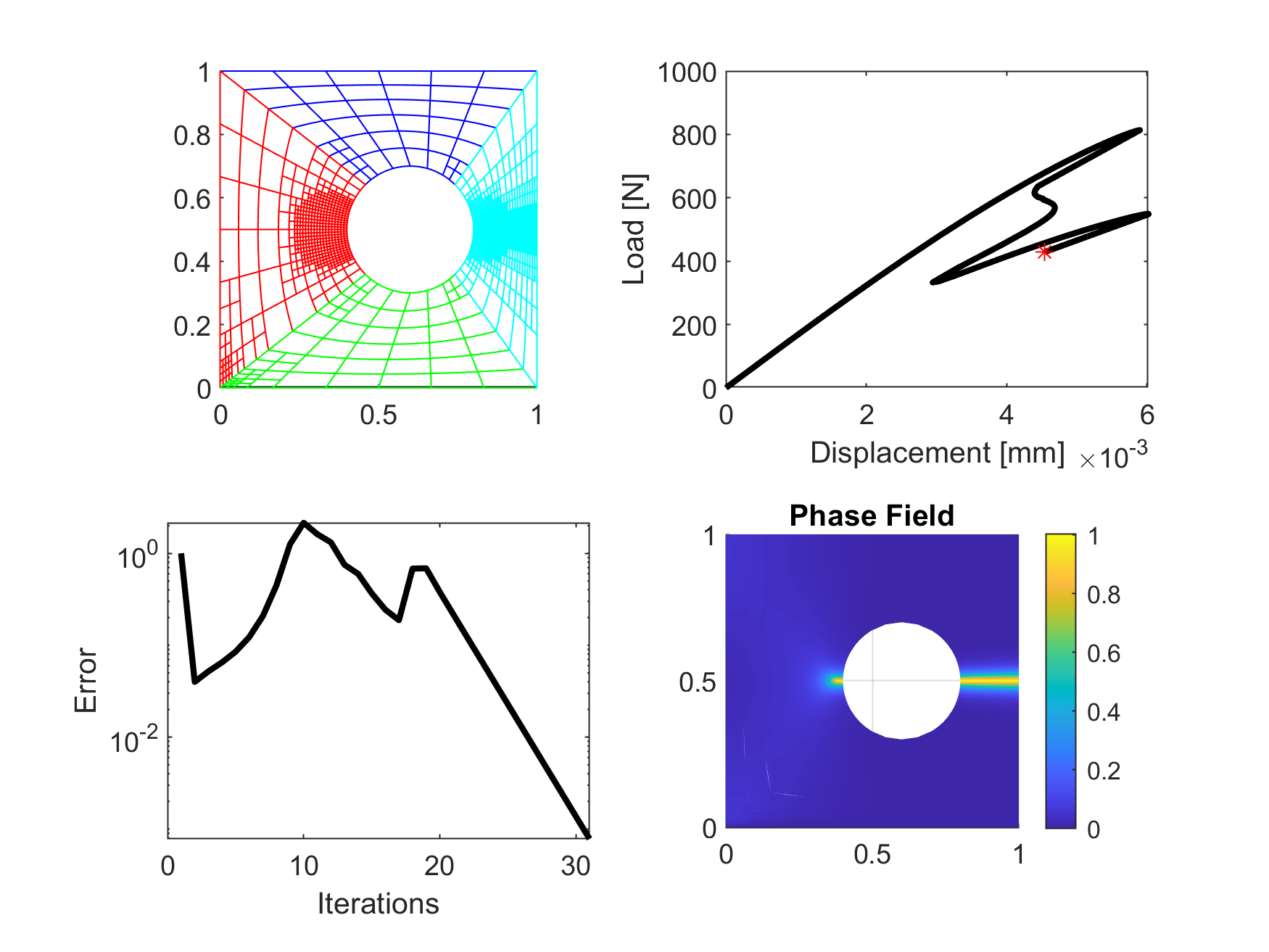}};
    \end{tikzpicture}
    \caption{ }
  \end{subfigure}
  \begin{subfigure}[t]{0.22\textwidth}
  \centering
    \begin{tikzpicture}[scale=0.5]
    \node[inner sep=0pt] () at (0,0)
    {\includegraphics[width=3.0cm,trim=8.35cm 1.5cm 2.65cm 6.25cm, clip]{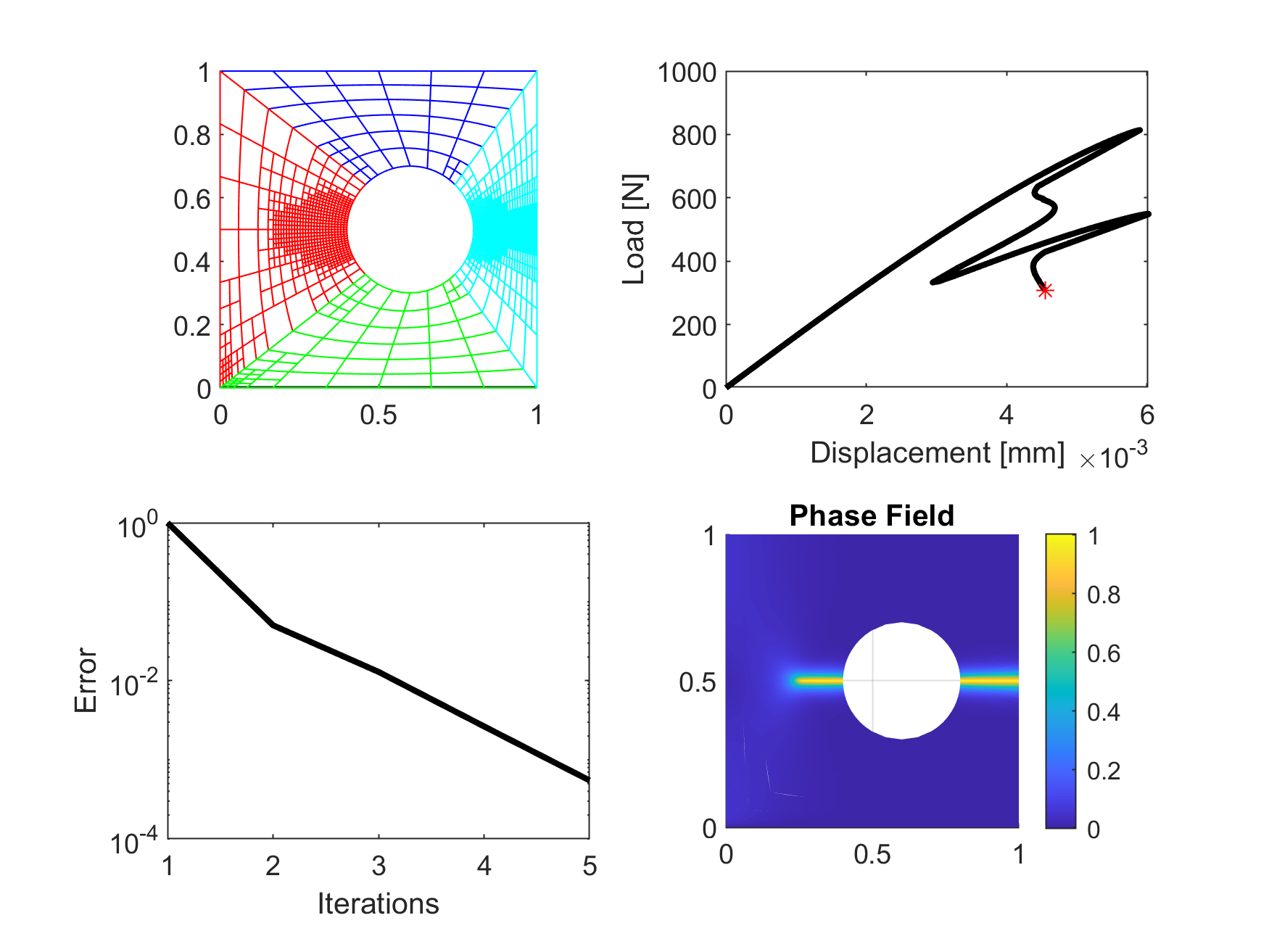}};
    \end{tikzpicture}
    \caption{ }
  \end{subfigure}
  \begin{subfigure}[t]{0.22\textwidth}
  \centering
    \begin{tikzpicture}[scale=0.5]
    \node[inner sep=0pt] () at (0,0)
    {\includegraphics[width=3.0cm,trim=8.35cm 1.5cm 2.65cm 6.25cm, clip]{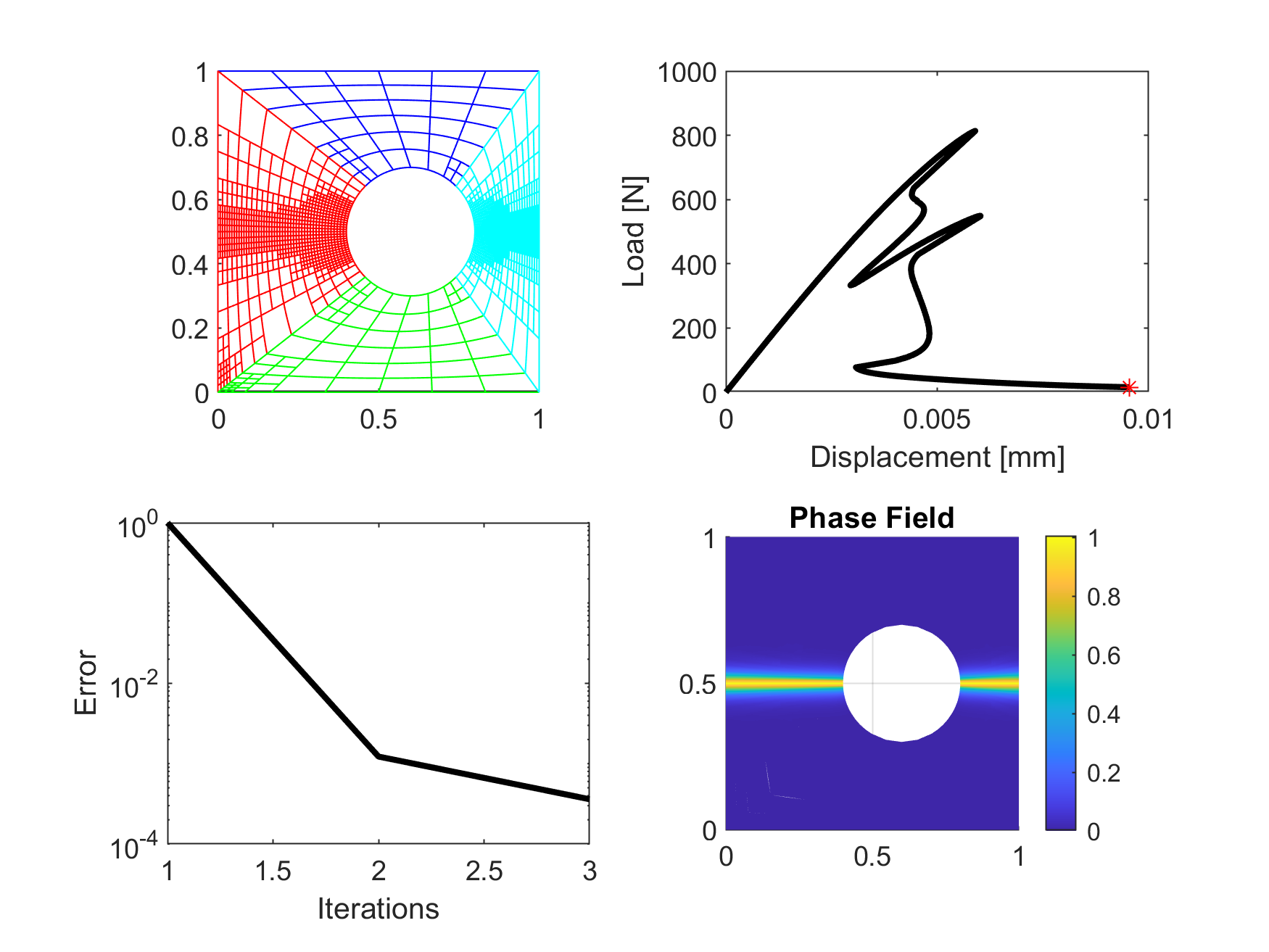}};
    \end{tikzpicture}
    \caption{ }
  \end{subfigure}
  \begin{subfigure}[t]{0.08\textwidth}
    \begin{tikzpicture}
    \begin{axis}[%
    hide axis,
    scale only axis,
    height=1.75\linewidth,
    width=0.01\linewidth,
    point meta min=0.0,
    point meta max=1.0,
    % colormap/bluered,                     % Colormap preset
    colorbar right,                  % Active colorbar
    colorbar sampled,                     % Steps in colorbar
    colorbar style={
        separate axis lines,
        samples=256,                        % Number of steps
    },
    every colorbar/.append style={
    height=\pgfkeysvalueof{/pgfplots/parent axis height},
    ytick={0,0.5,1},yticklabels={0,0.5\:\,$\pf$,1}
    }
    ]
    \addplot [draw=none] coordinates {(0,0)};
    \end{axis}
    \end{tikzpicture}
    \end{subfigure}
  % Mesh refinement plots    
  \begin{subfigure}[t]{0.23\textwidth}
  \centering
    \begin{tikzpicture}[scale=0.5]
    \node[inner sep=0pt] () at (0,0)
    {\includegraphics[width=3.0cm,trim=2.45cm 6.40cm 8.15cm 0.5cm, clip]{Figures/EH/Ref/Step39.png}};
    \end{tikzpicture}
    \caption{ }
  \end{subfigure}
  %\hfill
  % 2nd
  \begin{subfigure}[t]{0.22\textwidth}
  \centering
    \begin{tikzpicture}[scale=0.5]
    \node[inner sep=0pt] () at (0,0)
    {\includegraphics[width=3.0cm,trim=2.45cm 6.40cm 8.15cm 0.5cm, clip]{Figures/EH/Ref/Step97.png}};
    \end{tikzpicture}
    \caption{ }
  \end{subfigure}
  \begin{subfigure}[t]{0.22\textwidth}
  \centering
    \begin{tikzpicture}[scale=0.5]
    \node[inner sep=0pt] () at (0,0)
    {\includegraphics[width=3.0cm,trim=2.45cm 6.40cm 8.15cm 0.5cm, clip]{Figures/EH/Ref/Step127.png}};
    \end{tikzpicture}
    \caption{ }
  \end{subfigure}
  \begin{subfigure}[t]{0.22\textwidth}
  \centering
    \begin{tikzpicture}[scale=0.5]
    \node[inner sep=0pt] () at (0,0)
    {\includegraphics[width=3.0cm,trim=2.45cm 6.40cm 8.15cm 0.5cm, clip]{Figures/EH/Ref/Step209.png}};
    \end{tikzpicture}
    \caption{ }
  \end{subfigure}
  \caption{Figures (a-d) shows the evolution of the phase-field, and the corresponding refined meshes are shown in Figures (e-h) for the eccentric hole. The colors in the latter figures represent the different patches in the IGA context.}
  \label{sec5:fig:adap_mesh_EH}
\end{figure}

Figure \ref{sec5:fig:EH_lodi} presents the load-displacement curves for the unit square specimen with an eccentric hole under tension. The different curves correspond to the choice of the degradation function, quadratic and cubic. It is observed that the specimen exhibits a linear pre-peak behaviour with the cubic degradation function for $s<1$. This linear stage is missing for the quadratic degradation function and the cubic degradation function with $s=1$. However, irrespective of the choice of the degradation function, two snap-back behaviours are observed. The first one occurring at the onset of the localization from the hole towards the right of the specimen. The second snap-back behaviour corresponds to the onset of the localization from the hole towards the left edge of the specimen. Next, in Figure \ref{sec5:fig:EH_failure}, the phase-field fracture topology at the final step of the analysis is presented. The fracture topology similar to those observed in the literature \cite{May2016}. The refined meshes corresponding to the different stages in the evolution of the phase-field is shown in Figure \ref{sec5:fig:adap_mesh_EH}. 

% Figures (a-d) shows the evolution of the phase-field, and the corresponding refined meshes are shown in Figures (e-h) for the eccentric hole. The colors in the latter figures represent the different patches in the IGA context. 

\subsection{Single Edge Notched specimen under Tension (SENT)}\label{sec5:SENtension}

A unit square (in mm) embedded with a horizontal notch, midway along the height is considered, as shown in Figure \ref{sec5:fig:SENTdiagram}. The length of the notch is equal to half of the edge length of the plate (shown in red). The notch is modelled explicitly in the finite element mesh. A quasi-static loading is applied at the top boundary in the form of prescribed displacement increment $\tilde{\lambda}\hat{\bdisp}$, where $\hat{\bdisp}$ is a unit load vector and $\tilde{\lambda}$ is the load factor. When the analysis is started, the displacement-control approach is adopted and $\hat{\lambda}$ is incremented in steps of $1e-4$. Following the switch to the arc-length method, $\hat{\lambda}$ becomes an unknown variable, and is solved with the arc-length constraint equation. Furthermore, the bottom boundary remains fixed. The additional parameters required for the analysis are presented in Table \ref{sec5:table:SENTparams}. 

\begin{figure}[ht]
\begin{minipage}[b]{0.45\linewidth}
\centering
\begin{tikzpicture}[scale=0.7]
    \coordinate (K) at (0,0);
    % Square
    \draw[line width=0.75pt,black] (-2.5,-2.5) to (2.5,-2.5);
    \draw[line width=0.75pt,black] (2.5,-2.5) to (2.5,2.5);
    \draw[line width=0.75pt,black] (2.5,2.5) to (-2.5,2.5);
    \draw[line width=0.75pt,black] (-2.5,2.5) to (-2.5,-2.5);
    \draw[line width=1.5pt,red] (-2.5,0) to (0,0);
    % Top Boundary
    \draw[line width=0.75pt,black] (2.5,2.75) to (-2.5,2.75);
    \draw[->,line width=1.5pt,black] (0.0,3.0) to (0.0,3.5);
    \node[ ] at (0,3.85) {$\tilde{\lambda}\hat{\bdisp}$};
    % Bottom Boundary
    \draw[line width=0.75pt,black] (-2.5,-2.5) to (-2.75,-2.75);
    \draw[line width=0.75pt,black] (-2.25,-2.5) to (-2.5,-2.75);
    \draw[line width=0.75pt,black] (-2.0,-2.5) to (-2.25,-2.75);
    \draw[line width=0.75pt,black] (-1.75,-2.5) to (-2.0,-2.75);
    \draw[line width=0.75pt,black] (-1.5,-2.5) to (-1.75,-2.75);
    \draw[line width=0.75pt,black] (-1.25,-2.5) to (-1.5,-2.75);
    \draw[line width=0.75pt,black] (-1.0,-2.5) to (-1.25,-2.75);
    \draw[line width=0.75pt,black] (-0.75,-2.5) to (-1.0,-2.75);
    \draw[line width=0.75pt,black] (-0.5,-2.5) to (-0.75,-2.75);
    \draw[line width=0.75pt,black] (-0.25,-2.5) to (-0.5,-2.75);
    \draw[line width=0.75pt,black] (0.0,-2.5) to (-0.25,-2.75);
    \draw[line width=0.75pt,black] (0.25,-2.5) to (0.0,-2.75);
    \draw[line width=0.75pt,black] (0.5,-2.5) to (0.25,-2.75);
    \draw[line width=0.75pt,black] (0.75,-2.5) to (0.5,-2.75);
    \draw[line width=0.75pt,black] (1.0,-2.5) to (0.75,-2.75);
    \draw[line width=0.75pt,black] (1.25,-2.5) to (1.0,-2.75);
    \draw[line width=0.75pt,black] (1.5,-2.5) to (1.25,-2.75);
    \draw[line width=0.75pt,black] (1.75,-2.5) to (1.5,-2.75);
    \draw[line width=0.75pt,black] (2.0,-2.5) to (1.75,-2.75);
    \draw[line width=0.75pt,black] (2.25,-2.5) to (2.0,-2.75);
    \draw[line width=0.75pt,black] (2.5,-2.5) to (2.25,-2.75);
    % \node[ ] at (-0.5,-2.05) {$\tilde{\bdisp}_{p} = 0$};
    \end{tikzpicture}
\caption{SENT experiment}
\label{sec5:fig:SENTdiagram}
\end{minipage}
\begin{minipage}[b]{0.45\linewidth}
\centering
\begin{tabular}{ll} \hline
  \textbf{Parameters} & \textbf{Value} \\ \hline 
  Model & AT2 \\
  $\Psi^f$ & No Split \\
  $\lambda$ & 121.154 [GPa] \\
  $\mu$ & 80.769 [GPa] \\
  $\gc$ & 2700 [N/m] \\
  $\l$ & 2e-2 [mm] \\
  $\pf_{\text{threshold}}$ & 0.2 \\  
  $\Delta\tau_{max}$ & 0.025 [N] \\ \hline
  \end{tabular}
\captionof{table}{Model parameters}
\label{sec5:table:SENTparams}
\end{minipage}
\end{figure}

\begin{figure}[!ht]
  \begin{subfigure}[t]{0.45\textwidth}
  \centering
    \begin{tikzpicture}[thick,scale=0.95, every node/.style={scale=0.95}]
    \begin{axis}[legend style={draw=none}, legend columns = 2,
      transpose legend, ylabel={Load\:[N]},xlabel={Displacement\:[mm]}, xmin=0, ymin=0, xmax=0.0075, ymax=1300, yticklabel style={/pgf/number format/.cd,fixed,precision=2},
                 every axis plot/.append style={very thick}]
    \pgfplotstableread[col sep = comma]{./Data/SENT/lodi_quad.txt}\Adata;
    \pgfplotstableread[col sep = comma]{./Data/SENT/lodi_cubics0_01.txt}\Bdata;
    \pgfplotstableread[col sep = comma]{./Data/SENT/lodi_cubics0_1.txt}\Cdata;
    \pgfplotstableread[col sep = comma]{./Data/SENT/lodi_cubics1.txt}\Ddata;
    \addplot [ 
           color=black, 
%           only marks, 
           mark=*, 
           mark size=0.25pt, 
         ]
         table
         [
           x expr=\thisrowno{1}, 
           y expr=\thisrowno{0}
         ] {\Adata};
         \addlegendentry{Quadratic}
%    \addplot [ 
%           color=red, 
%           only marks, 
%           mark=*, 
%           mark size=0.25pt, 
%         ]
%         table
%         [
%           x expr=\thisrowno{1}, 
%           y expr=\thisrowno{0}
%         ] {\Bdata};
%         \addlegendentry{Cubic: $s = 0.01$}     
    \addplot [ 
           color=red, 
%           only marks, 
           mark=*, 
           mark size=0.25pt, 
         ]
         table
         [
           x expr=\thisrowno{1}, 
           y expr=\thisrowno{0}
         ] {\Cdata};
         \addlegendentry{Cubic: $s = 0.1$}     
    \addplot [ 
           color=blue, 
%           only marks, 
           mark=*, 
           mark size=0.25pt, 
         ]
         table
         [
           x expr=\thisrowno{1}, 
           y expr=\thisrowno{0}
         ] {\Ddata};
         \addlegendentry{Cubic: $s = 1.0$}          
    \end{axis}
    \end{tikzpicture}
    \caption{ }
    \label{sec5:fig:SENT_lodi}
  \end{subfigure}
  \hfill
  \begin{subfigure}[t]{0.45\textwidth}
  \centering
    \begin{tikzpicture}
    \node[inner sep=0pt] () at (0,0)
    {\includegraphics[width=5.5cm,trim=8.35cm 1.25cm 2.65cm 6.25cm, clip]{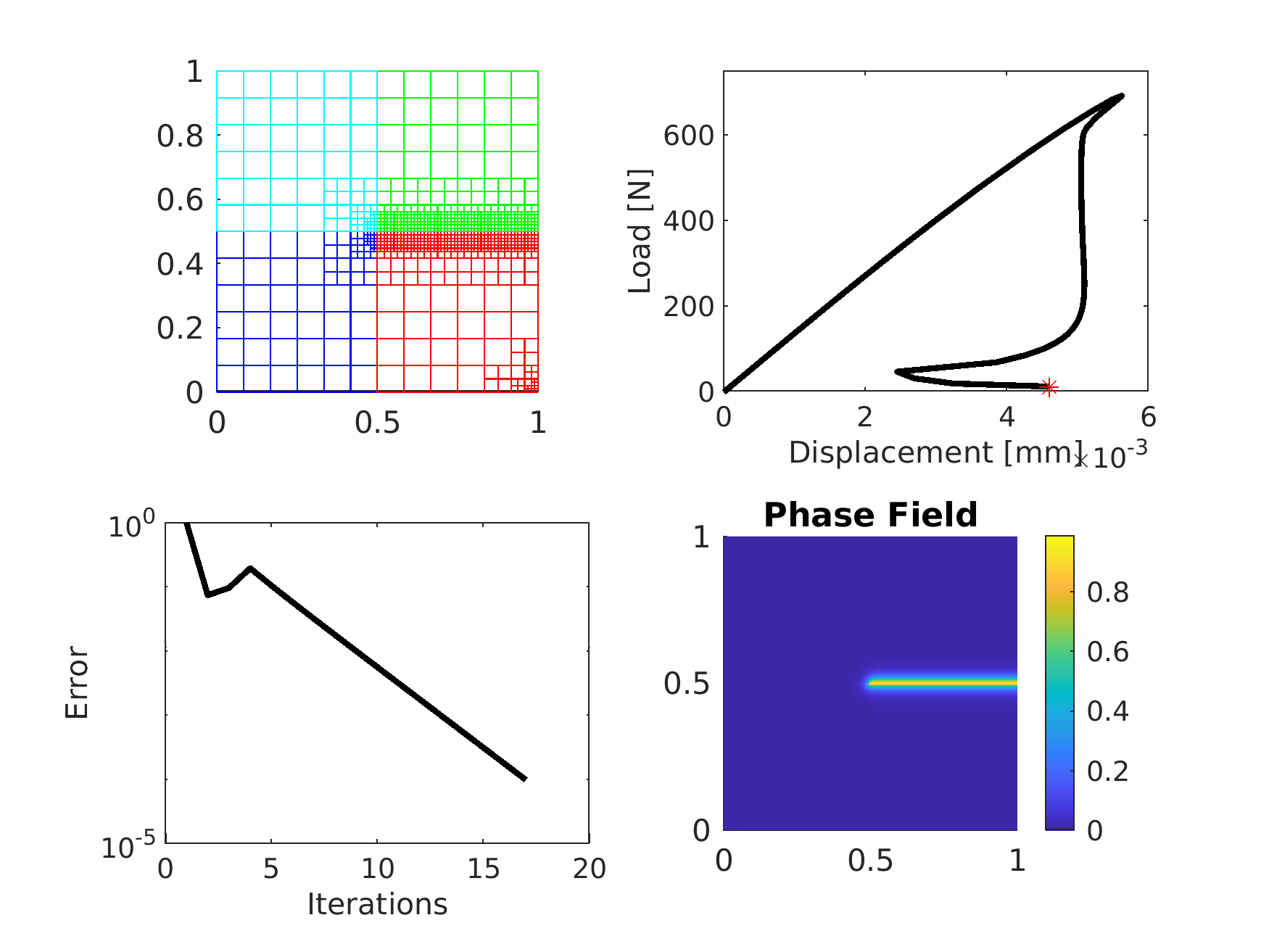}};
    \node[inner sep=0pt] () at (-1.3,-3.05)
    {\begin{axis}[
    hide axis,
    scale only axis,
    height=0pt,
    width=0pt,
    %colormap/bluered,
    colorbar horizontal,
    point meta min=0,
    point meta max=1,
    colorbar style={
        width=4.85cm,
        xtick={0,0.5,1.0},
        xticklabel style = {yshift=-0.075cm}
    }]
    \addplot [draw=none] coordinates {(0,0)};
    \end{axis}};
    \node[inner sep=0pt] () at (0,-3.75) {$\pf$};
    \end{tikzpicture}
    \caption{ }
    \label{sec5:fig:SENT_failure}
  \end{subfigure}
  \caption{Figure (a) presents the load-displacement curves for the single edge notched specimen under tension. The legend entries correspond to the choice of degradation functions. Figure (b) shows the distribution of the phase-field variable at the final step of the analysis.}
\end{figure}

\begin{figure}[!ht]
 \begin{subfigure}[t]{0.22\textwidth}
  \centering
    \begin{tikzpicture}[scale=0.5]
    \node[inner sep=0pt] () at (0,0)
    {\includegraphics[width=3.0cm,trim=8.35cm 1.5cm 2.65cm 6.25cm, clip]{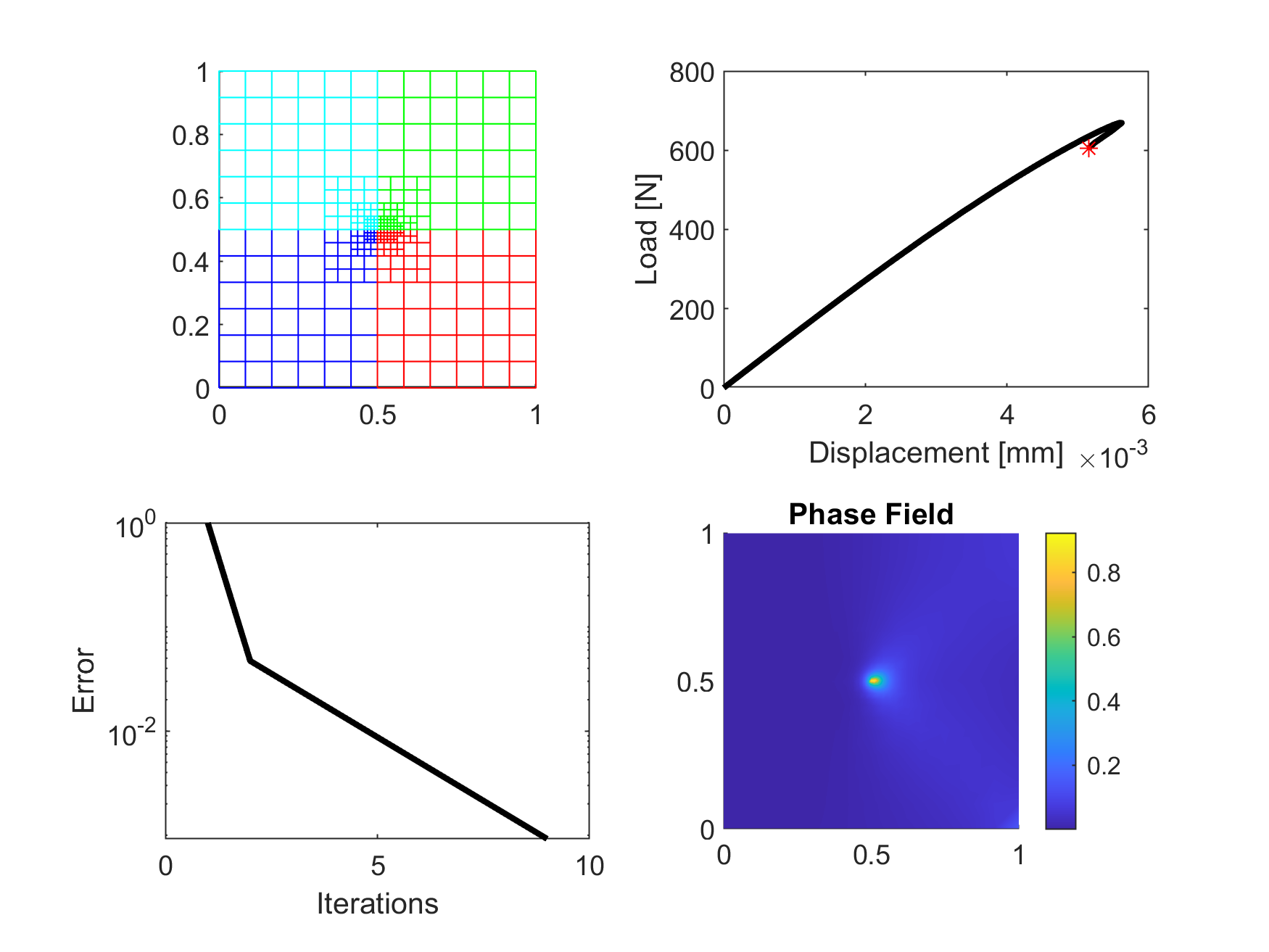}};
    \end{tikzpicture}
    \caption{ }
  \end{subfigure}
  %\hfill
  %
  \begin{subfigure}[t]{0.22\textwidth}
  \centering
    \begin{tikzpicture}[scale=0.5]
    \node[inner sep=0pt] () at (0,0)
    {\includegraphics[width=3.0cm,trim=8.35cm 1.5cm 2.65cm 6.25cm, clip]{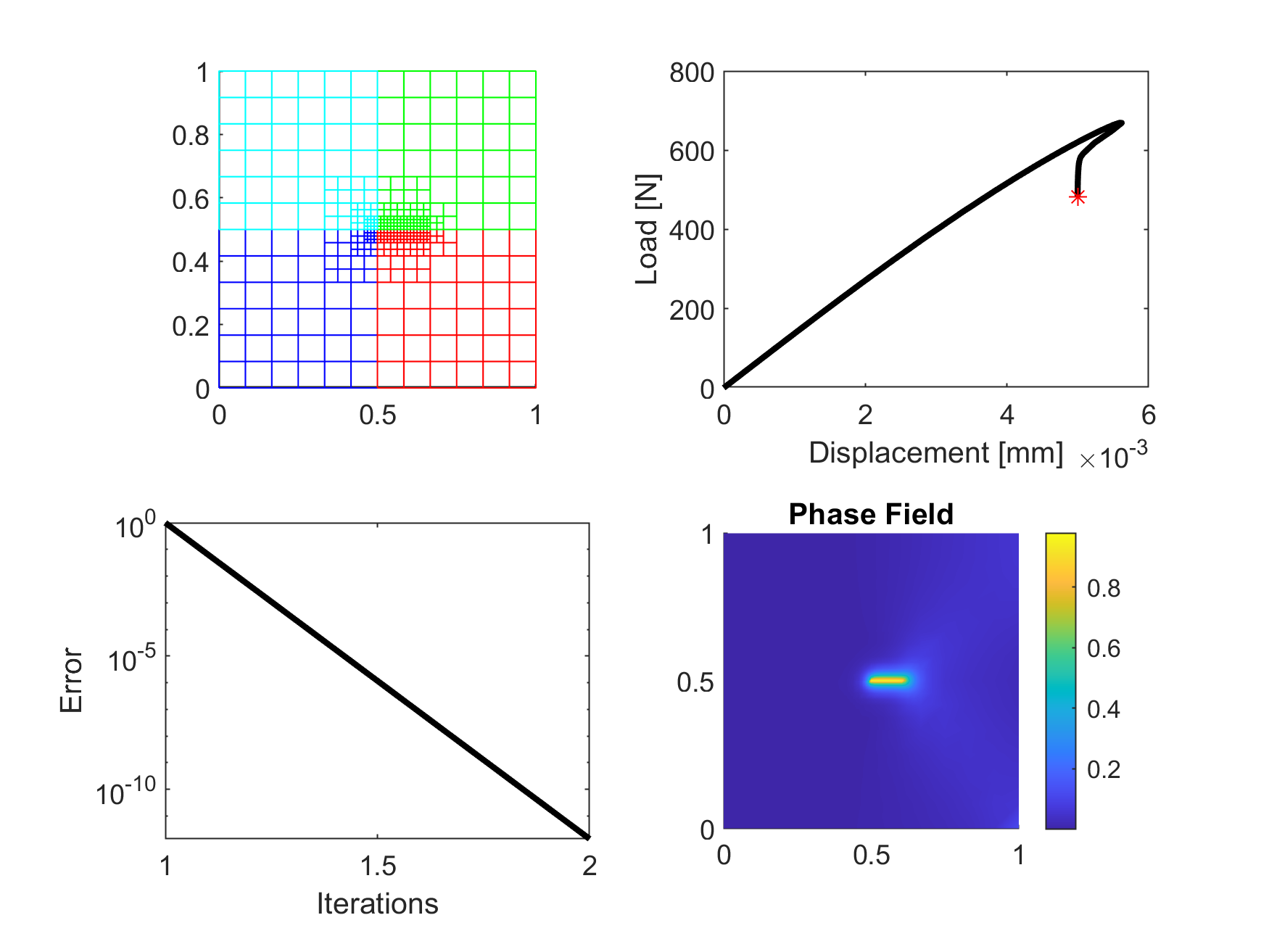}};
    \end{tikzpicture}
    \caption{ }
  \end{subfigure}
  \begin{subfigure}[t]{0.22\textwidth}
  \centering
    \begin{tikzpicture}[scale=0.5]
    \node[inner sep=0pt] () at (0,0)
    {\includegraphics[width=3.0cm,trim=8.35cm 1.5cm 2.65cm 6.25cm, clip]{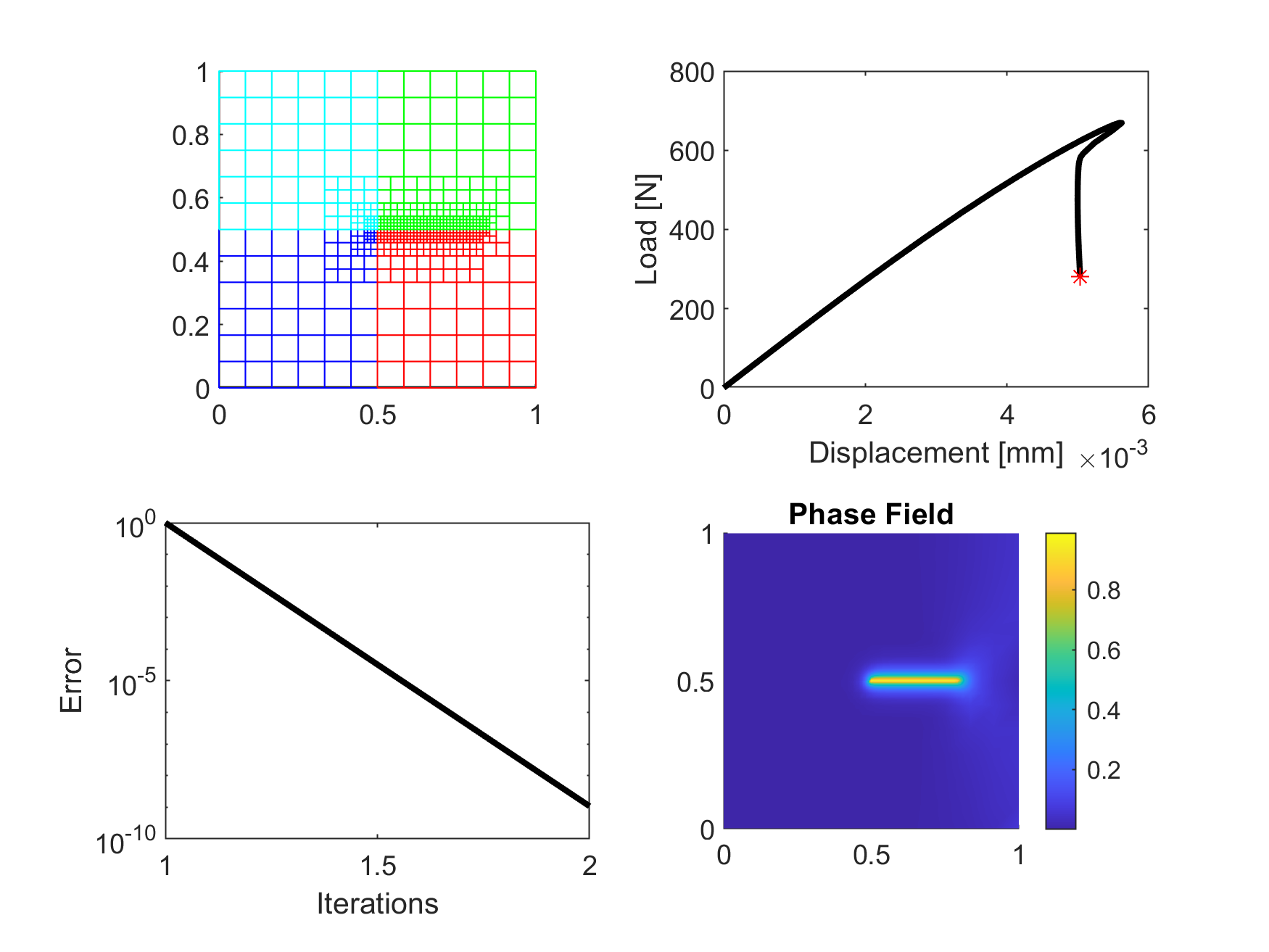}};
    \end{tikzpicture}
    \caption{ }
  \end{subfigure}
  \begin{subfigure}[t]{0.22\textwidth}
  \centering
    \begin{tikzpicture}[scale=0.5]
    \node[inner sep=0pt] () at (0,0)
    {\includegraphics[width=3.0cm,trim=8.35cm 1.5cm 2.65cm 6.25cm, clip]{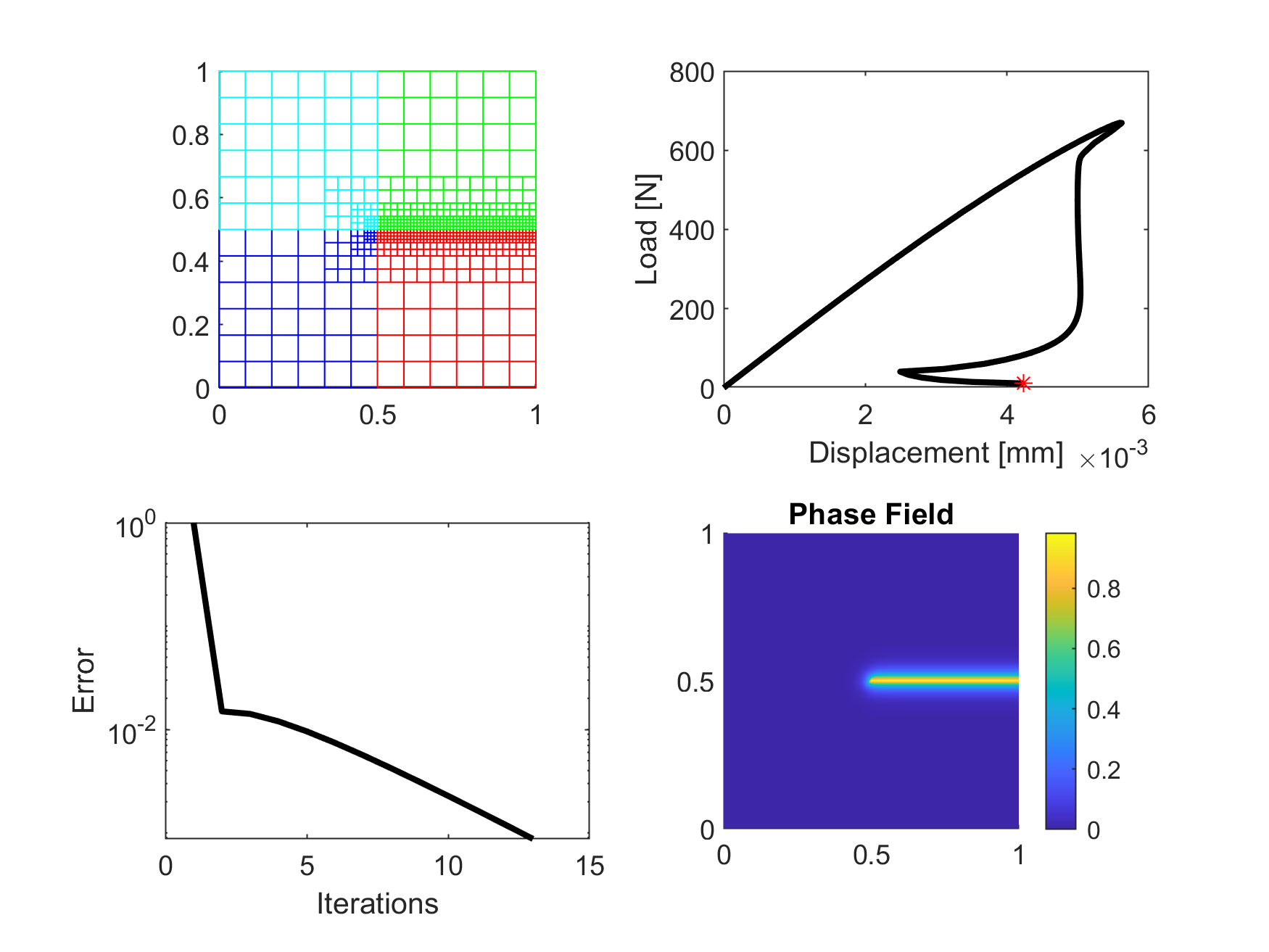}};
    \end{tikzpicture}
    \caption{ }
  \end{subfigure}
  \begin{subfigure}[t]{0.08\textwidth}
    \begin{tikzpicture}
    \begin{axis}[%
    hide axis,
    scale only axis,
    height=1.75\linewidth,
    width=0.01\linewidth,
    point meta min=0.0,
    point meta max=1.0,
    % colormap/bluered,                     % Colormap preset
    colorbar right,                  % Active colorbar
    colorbar sampled,                     % Steps in colorbar
    colorbar style={
        separate axis lines,
        samples=256,                        % Number of steps
    },
    every colorbar/.append style={
    height=\pgfkeysvalueof{/pgfplots/parent axis height},
    ytick={0,0.5,1},yticklabels={0,0.5\:\,$\pf$,1}
    }
    ]
    \addplot [draw=none] coordinates {(0,0)};
    \end{axis}
    \end{tikzpicture}
    \end{subfigure}
  % Mesh refinement plots    
  \begin{subfigure}[t]{0.23\textwidth}
  \centering
    \begin{tikzpicture}[scale=0.5]
    \node[inner sep=0pt] () at (0,0)
    {\includegraphics[width=3.0cm,trim=2.45cm 6.40cm 8.15cm 0.5cm, clip]{Figures/SENT/Ref/Step32.png}};
    \end{tikzpicture}
    \caption{ }
  \end{subfigure}
  %\hfill
  % 2nd
  \begin{subfigure}[t]{0.22\textwidth}
  \centering
    \begin{tikzpicture}[scale=0.5]
    \node[inner sep=0pt] () at (0,0)
    {\includegraphics[width=3.0cm,trim=2.45cm 6.40cm 8.15cm 0.5cm, clip]{Figures/SENT/Ref/Step55.png}};
    \end{tikzpicture}
    \caption{}
  \end{subfigure}
  \begin{subfigure}[t]{0.22\textwidth}
  \centering
    \begin{tikzpicture}[scale=0.5]
    \node[inner sep=0pt] () at (0,0)
    {\includegraphics[width=3.0cm,trim=2.45cm 6.40cm 8.15cm 0.5cm, clip]{Figures/SENT/Ref/Step100.png}};
    \end{tikzpicture}
    \caption{ }
  \end{subfigure}
  \begin{subfigure}[t]{0.22\textwidth}
  \centering
    \begin{tikzpicture}[scale=0.5]
    \node[inner sep=0pt] () at (0,0)
    {\includegraphics[width=3.0cm,trim=2.45cm 6.40cm 8.15cm 0.5cm, clip]{Figures/SENT/Ref/Step148.png}};
    \end{tikzpicture}
    \caption{ }
  \end{subfigure}
  \caption{Figures (a-d) shows the evolution of the phase-field, and the corresponding refined meshes are shown in Figures (e-h) for the single edge notched specimen under tension. The colors in the latter figures represent the different patches in the IGA context.}
  \label{sec5:fig:adap_mesh_SENT}
\end{figure}

Figure \ref{sec5:fig:SENT_lodi} presents the load-displacement curves for the single edge notched specimen under tension. The different curves correspond to the choice of the degradation function, quadratic and cubic. Similar to the previous section, it is observed that the specimen exhibits a linear pre-peak behaviour with the cubic degradation function for $s<1$. The quadratic degradation function does not exhibit a linear pre-peak behaviour. Moreover, beyond the first snap-back behaviour, the post-peak branches of all curves are similar. Next, in Figure \ref{sec5:fig:SENT_failure}, the phase-field fracture topology at the final step of the analysis is presented. The fracture topology is consistent with those observed in the literature, for instance, \cite{Miehe2010}. The refined meshes corresponding to the different stages in the evolution of the phase-field is shown in Figure \ref{sec5:fig:adap_mesh_SENT}. 

Furthermore, sensitivity studies are carried out w.r.t to the choice of the prescribed maximum dissipation $\Delta\tau_{max}$ and the phase-field threshold for mesh refinement. The results in Appendix \ref{appA:sensitivity} present the $\Delta\tau_{max}$ for which similar load-displacement curves are obtained. Next, it is also observed that the adaptive mesh refinement technique with different values of $\pf_{th.}$ yield similar load-displacement curve as that obtained on a fixed mesh. Finally, in Appendix \ref{appA:comparison}, the proposed monolithic solver is compared with conventional alternate minimization solver and the quasi-Newton Raphson method \cite{Heister2015}. For the SENT problem, it is observed that all methods/solvers yield similar peak load and pre-peak behaviour. However, the proposed monolithic solver augmented with the arc-length method captures snap-back behaviours, which is not possible with alternative minimization and quasi-Newton Raphson method.

\subsection{Single Edge Notched specimen under Shear (SENS)}

The shear test is carried out on the single edge notched specimen by loading in the horizontal direction, as shown in Figure \ref{sec5:fig:SENSdiagram}. The model parameters remain same as presented in Table \ref{sec5:table:SENSparams}. Similar to the SENT model, a quasi-static loading is applied at the top boundary in the form of prescribed displacement increment $\tilde{\lambda}\hat{\bdisp}$, where $\hat{\bdisp}$ is a unit load vector and $\tilde{\lambda}$ is the load factor. The loading in applied in the form of prescribed displacement increment $\tilde{\lambda}\hat{\bdisp}$, where $\hat{\bdisp}$ is a unit load vector and $\tilde{\lambda}$ is the load factor. When the analysis is started, the displacement-control approach is adopted and $\hat{\lambda}$ is incremented in steps of $5e-4$ [mm]. Following the switch to the arc-length method, $\hat{\lambda}$ becomes an unknown variable, and is solved with the arc-length constraint equation. The bottom boundary of the specimen is fixed, while the left and the right edges are restricted in the vertical direction.

\begin{figure}[ht]
\begin{minipage}[b]{0.45\linewidth}
\centering
\begin{tikzpicture}[scale=0.7]
    \coordinate (K) at (0,0);
    % Square
    \draw[line width=0.75pt,black] (-2.5,-2.5) to (2.5,-2.5);
    \draw[line width=0.75pt,black] (2.5,-2.5) to (2.5,2.5);
    \draw[line width=0.75pt,black] (2.5,2.5) to (-2.5,2.5);
    \draw[line width=0.75pt,black] (-2.5,2.5) to (-2.5,-2.5);
    \draw[line width=1.5pt,red] (-2.5,0) to (0,0);
    % Top Boundary
    \draw[line width=0.75pt,black] (2.5,2.75) to (-2.5,2.75);
    \draw[->,line width=1.5pt,black] (0.1,3.25) to (0.8,3.25);
    \node[ ] at (-0.5,3.25) {$\tilde{\lambda}\hat{\bdisp}$};
    %\node[ ] at (0,3.85) {$\disp_y=0$};
    % Bottom Boundary
    \draw[line width=0.75pt,black] (-2.5,-2.5) to (-2.75,-2.75);
    \draw[line width=0.75pt,black] (-2.25,-2.5) to (-2.5,-2.75);
    \draw[line width=0.75pt,black] (-2.0,-2.5) to (-2.25,-2.75);
    \draw[line width=0.75pt,black] (-1.75,-2.5) to (-2.0,-2.75);
    \draw[line width=0.75pt,black] (-1.5,-2.5) to (-1.75,-2.75);
    \draw[line width=0.75pt,black] (-1.25,-2.5) to (-1.5,-2.75);
    \draw[line width=0.75pt,black] (-1.0,-2.5) to (-1.25,-2.75);
    \draw[line width=0.75pt,black] (-0.75,-2.5) to (-1.0,-2.75);
    \draw[line width=0.75pt,black] (-0.5,-2.5) to (-0.75,-2.75);
    \draw[line width=0.75pt,black] (-0.25,-2.5) to (-0.5,-2.75);
    \draw[line width=0.75pt,black] (0.0,-2.5) to (-0.25,-2.75);
    \draw[line width=0.75pt,black] (0.25,-2.5) to (0.0,-2.75);
    \draw[line width=0.75pt,black] (0.5,-2.5) to (0.25,-2.75);
    \draw[line width=0.75pt,black] (0.75,-2.5) to (0.5,-2.75);
    \draw[line width=0.75pt,black] (1.0,-2.5) to (0.75,-2.75);
    \draw[line width=0.75pt,black] (1.25,-2.5) to (1.0,-2.75);
    \draw[line width=0.75pt,black] (1.5,-2.5) to (1.25,-2.75);
    \draw[line width=0.75pt,black] (1.75,-2.5) to (1.5,-2.75);
    \draw[line width=0.75pt,black] (2.0,-2.5) to (1.75,-2.75);
    \draw[line width=0.75pt,black] (2.25,-2.5) to (2.0,-2.75);
    \draw[line width=0.75pt,black] (2.5,-2.5) to (2.25,-2.75);
    % \node[ ] at (-0.5,-2.) {$\tilde{\bdisp}_{p} = 0$};
    % Left edge 1
    \draw[fill=gray!50] (-2.75,1.25) -- (-2.25,1.25) -- (-2.5,1.75)-- (-2.75,1.25);
    \draw[fill=black!75] (-2.5,1.05) circle (0.2);
    \draw[line width=1pt,black] (-2.25,0.8) to (-2.75,0.8);
    % Left edge 2
    \draw[fill=gray!50] (-2.75,-1.5) -- (-2.25,-1.5) -- (-2.5,-1.)-- (-2.75,-1.5);
    \draw[fill=black!75] (-2.5,-1.7) circle (0.2);
    \draw[line width=1pt,black] (-2.25,-1.95) to (-2.75,-1.95);
    % Right edge 1
    \draw[fill=gray!50] (2.75,1.25) -- (2.25,1.25) -- (2.5,1.75)-- (2.75,1.25);
    \draw[fill=black!75] (2.5,1.05) circle (0.2);
    \draw[line width=1pt,black] (2.25,0.8) to (2.75,0.8);
    % Right edge 2
    \draw[fill=gray!50] (2.75,-1.5) -- (2.25,-1.5) -- (2.5,-1.)-- (2.75,-1.5);
    \draw[fill=black!75] (2.5,-1.7) circle (0.2);
    \draw[line width=1pt,black] (2.25,-1.95) to (2.75,-1.95);
    \end{tikzpicture}
\caption{SENS experiment}
\label{sec5:fig:SENSdiagram}
\end{minipage}
\begin{minipage}[b]{0.45\linewidth}
\centering
\begin{tabular}{ll} \hline
  \textbf{Parameters} & \textbf{Value} \\ \hline 
  Model & AT2 \\
  $\Psi^f$ & Rankine \\
  $\lambda$ & 121.154 [GPa] \\
  $\mu$ & 80.769 [GPa] \\
  $\gc$ & 2700 [N/m] \\
  $\l$ & 2e-2 [mm] \\
  $\pf_{\text{threshold}}$ & 0.1 \\  
  $\Delta\tau_{max}$ & 0.025 [N] \\ \hline
  \end{tabular}
\captionof{table}{Model parameters}
\label{sec5:table:SENSparams}
\end{minipage}
\end{figure}

\begin{figure}[!ht]
 \begin{subfigure}[t]{0.45\textwidth}
  \centering
    \begin{tikzpicture}[thick,scale=0.95, every node/.style={scale=0.95}]
    \begin{axis}[legend style={draw=none}, legend columns = 2,
      transpose legend, ylabel={Load\:[N]},xlabel={Displacement\:[mm]}, xmin=0, ymin=0, xmax=0.014, ymax=800, yticklabel style={/pgf/number format/.cd,fixed,precision=2},
                 every axis plot/.append style={very thick}]
    \pgfplotstableread[col sep = comma]{./Data/SENS/lodi_quadn.txt}\Adata;
    %\pgfplotstableread[col sep = comma]{./Data/SENT/lodi_cubics0_01.txt}\Bdata;
    \pgfplotstableread[col sep = comma]{./Data/SENS/lodi_cubics0_1n.txt}\Cdata;
    \pgfplotstableread[col sep = comma]{./Data/SENS/lodi_cubics1n.txt}\Ddata;
    \addplot [ 
           color=black, 
%           only marks, 
           mark=*, 
           mark size=0.25pt, 
         ]
         table
         [
           x expr=\thisrowno{1}, 
           y expr=\thisrowno{0}
         ] {\Adata};
         \addlegendentry{Quadratic}
%    \addplot [ 
%           color=red, 
%           only marks, 
%           mark=*, 
%           mark size=0.25pt, 
%         ]
%         table
%         [
%           x expr=\thisrowno{1}, 
%           y expr=\thisrowno{0}
%         ] {\Bdata};
%         \addlegendentry{Cubic: $s = 0.01$}     
    \addplot [ 
           color=red, 
%           only marks, 
           mark=*, 
           mark size=0.25pt, 
         ]
         table
         [
           x expr=\thisrowno{1}, 
           y expr=\thisrowno{0}
         ] {\Cdata};
         \addlegendentry{Cubic: $s = 0.1$}     
    \addplot [ 
           color=blue, 
%           only marks, 
           mark=*, 
           mark size=0.25pt, 
         ]
         table
         [
           x expr=\thisrowno{1}, 
           y expr=\thisrowno{0}
         ] {\Ddata};
         \addlegendentry{Cubic: $s = 1.0$}          
    \end{axis}
    \end{tikzpicture}
    \caption{ }
    \label{sec5:fig:SENS_lodi}
  \end{subfigure}
  \hfill
  \begin{subfigure}[t]{0.45\textwidth}
  \centering
    \begin{tikzpicture}
    \node[inner sep=0pt] () at (0,0)
    {\includegraphics[width=5.5cm,trim=8.35cm 1.25cm 2.65cm 6.25cm, clip]{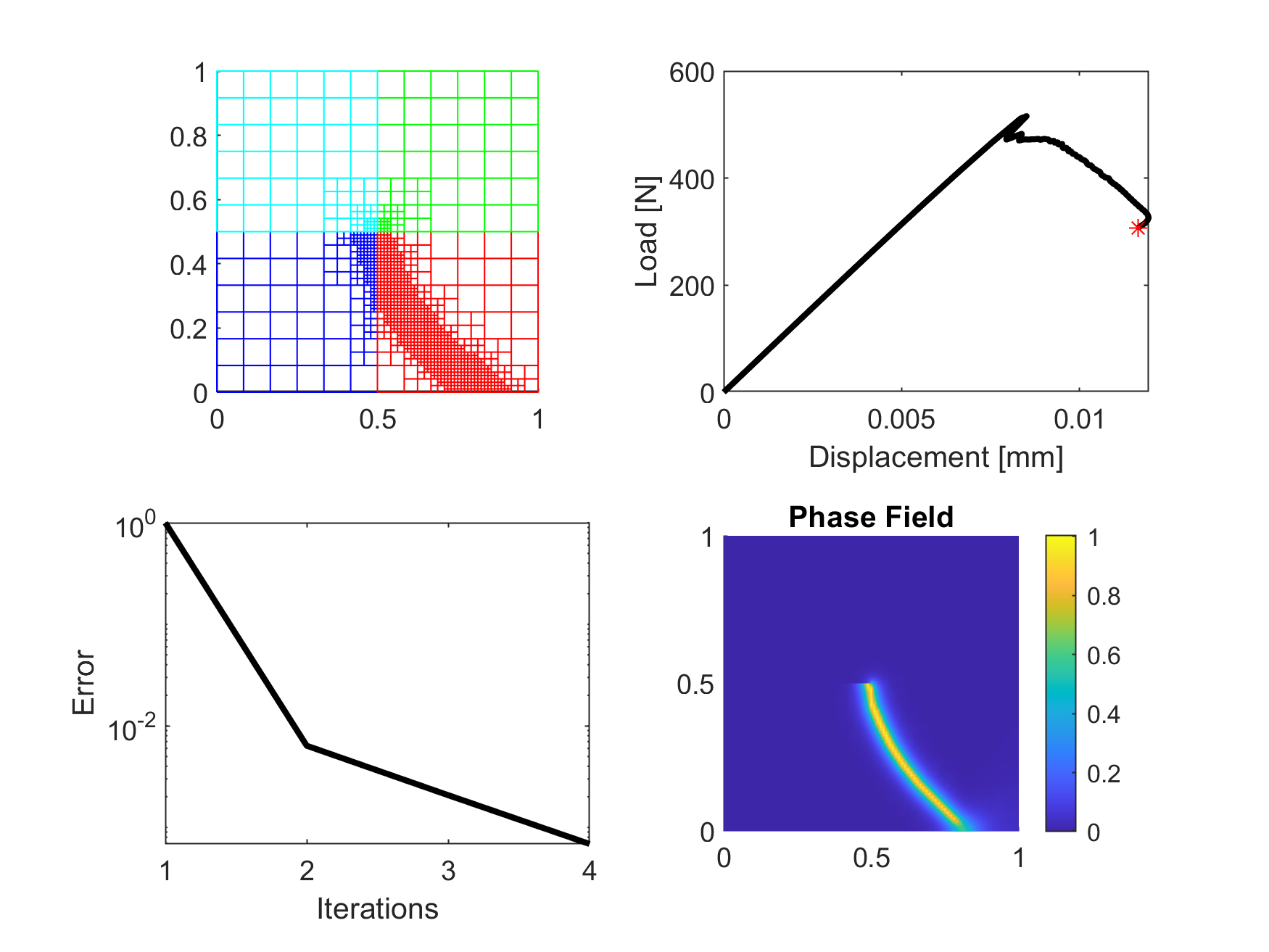}}; % Step261
    \node[inner sep=0pt] () at (-1.3,-3.05)
    {\begin{axis}[
    hide axis,
    scale only axis,
    height=0pt,
    width=0pt,
    %colormap/bluered,
    colorbar horizontal,
    point meta min=0,
    point meta max=1,
    colorbar style={
        width=4.85cm,
        xtick={0,0.5,1.0},
        xticklabel style = {yshift=-0.075cm}
    }]
    \addplot [draw=none] coordinates {(0,0)};
    \end{axis}};
    \node[inner sep=0pt] () at (0,-3.75) {$\pf$};
    \end{tikzpicture}
    \caption{ }
    \label{sec5:fig:SENS_failure}
  \end{subfigure}
  \caption{Figure (a) presents the load-displacement curves for the single edge notched specimen under shear. The legend entries correspond to the choice of degradation functions. Figure (b) shows the distribution of the phase-field variable at the final step of the analysis.}
\end{figure}

\begin{figure}[!ht]
 \begin{subfigure}[t]{0.22\textwidth}
  \centering
    \begin{tikzpicture}[scale=0.5]
    \node[inner sep=0pt] () at (0,0)
    {\includegraphics[width=3.0cm,trim=8.35cm 1.5cm 2.65cm 6.25cm, clip]{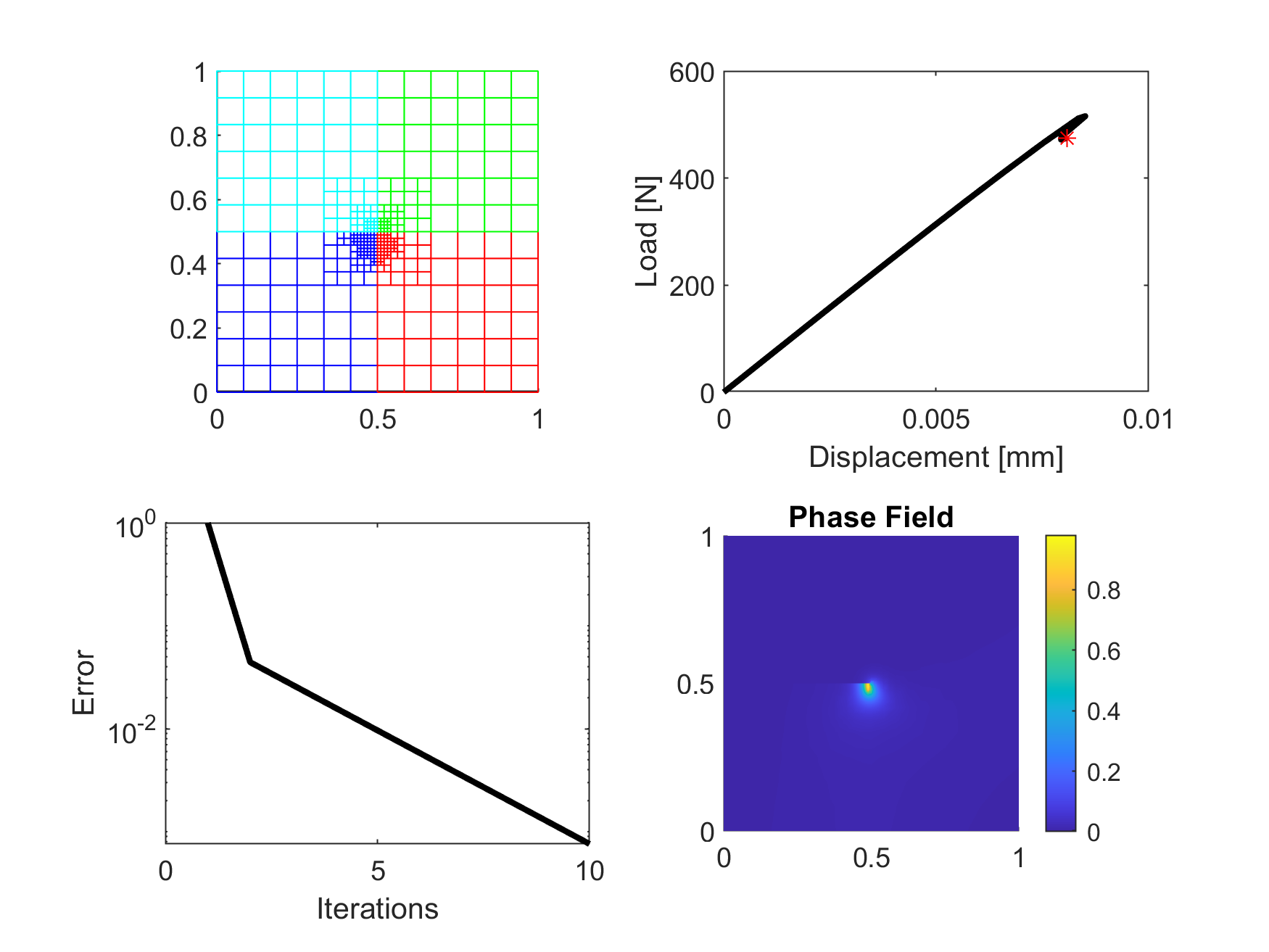}};
    \end{tikzpicture}
    \caption{ }
  \end{subfigure}
  %\hfill
  %
  \begin{subfigure}[t]{0.22\textwidth}
  \centering
    \begin{tikzpicture}[scale=0.5]
    \node[inner sep=0pt] () at (0,0)
    {\includegraphics[width=3.0cm,trim=8.35cm 1.5cm 2.65cm 6.25cm, clip]{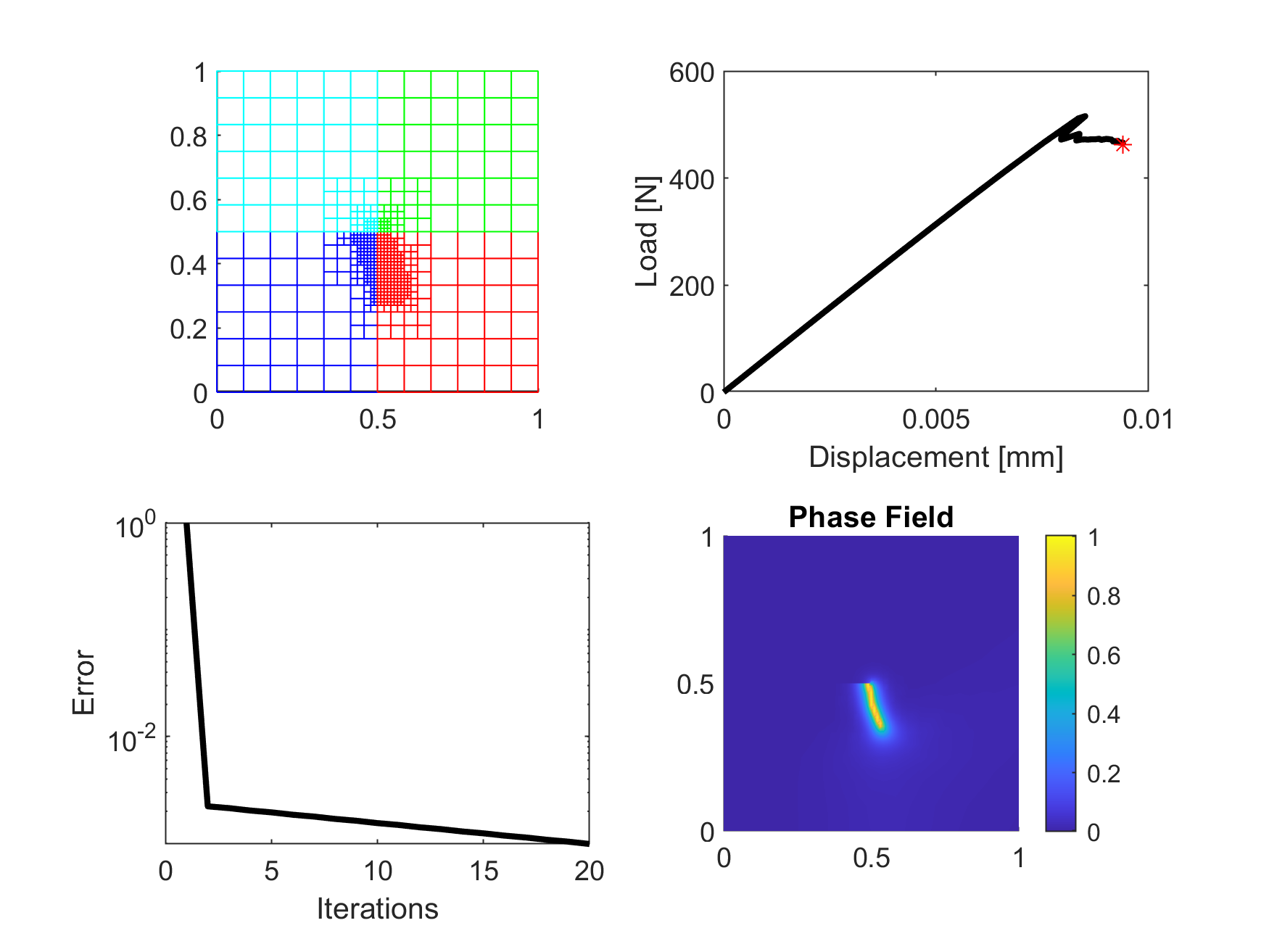}};
    \end{tikzpicture}
    \caption{ }
  \end{subfigure}
  \begin{subfigure}[t]{0.22\textwidth}
  \centering
    \begin{tikzpicture}[scale=0.5]
    \node[inner sep=0pt] () at (0,0)
    {\includegraphics[width=3.0cm,trim=8.35cm 1.5cm 2.65cm 6.25cm, clip]{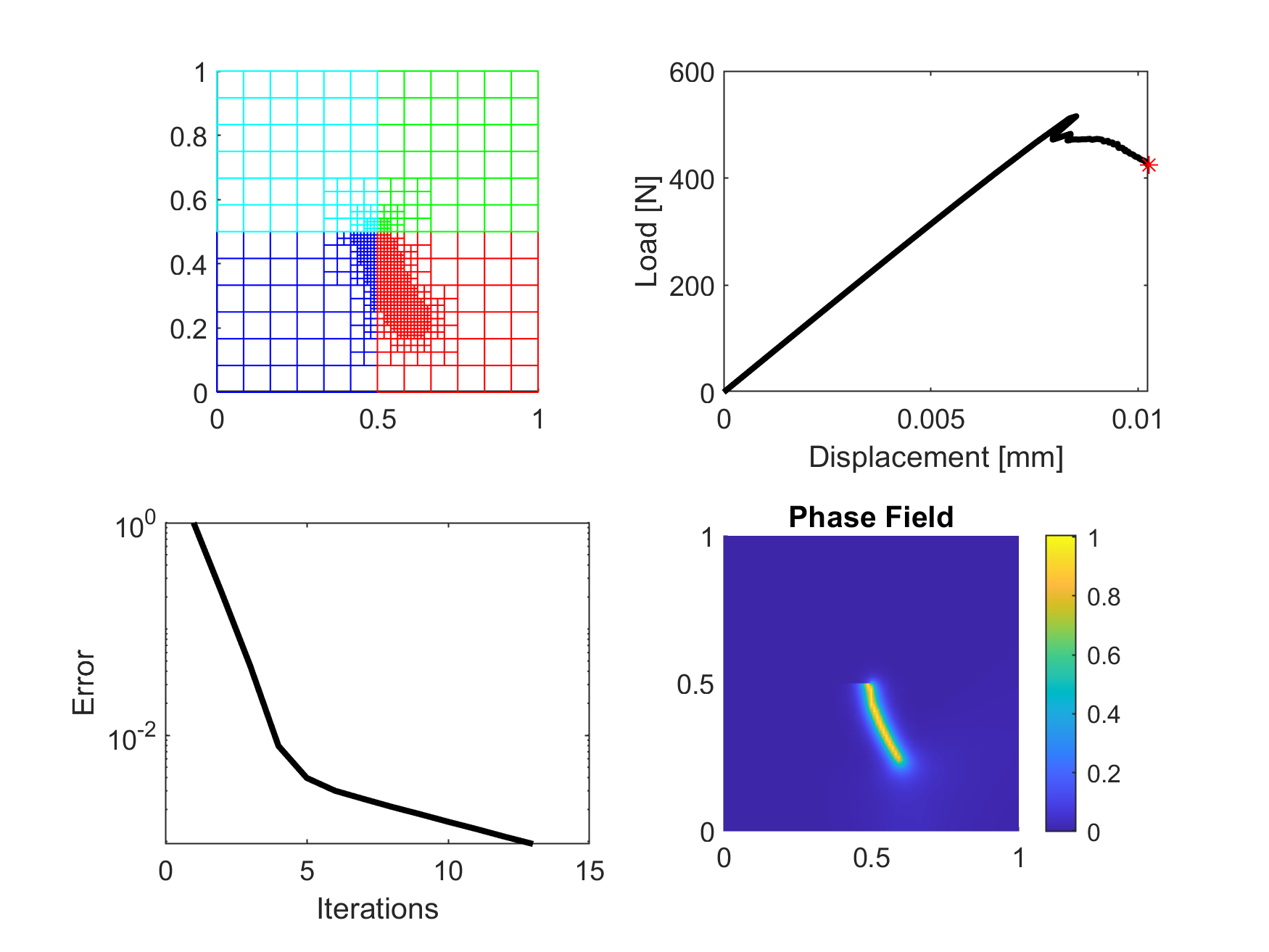}};
    \end{tikzpicture}
    \caption{ }
  \end{subfigure}
  \begin{subfigure}[t]{0.22\textwidth}
  \centering
    \begin{tikzpicture}[scale=0.5]
    \node[inner sep=0pt] () at (0,0)
    {\includegraphics[width=3.0cm,trim=8.35cm 1.5cm 2.65cm 6.25cm, clip]{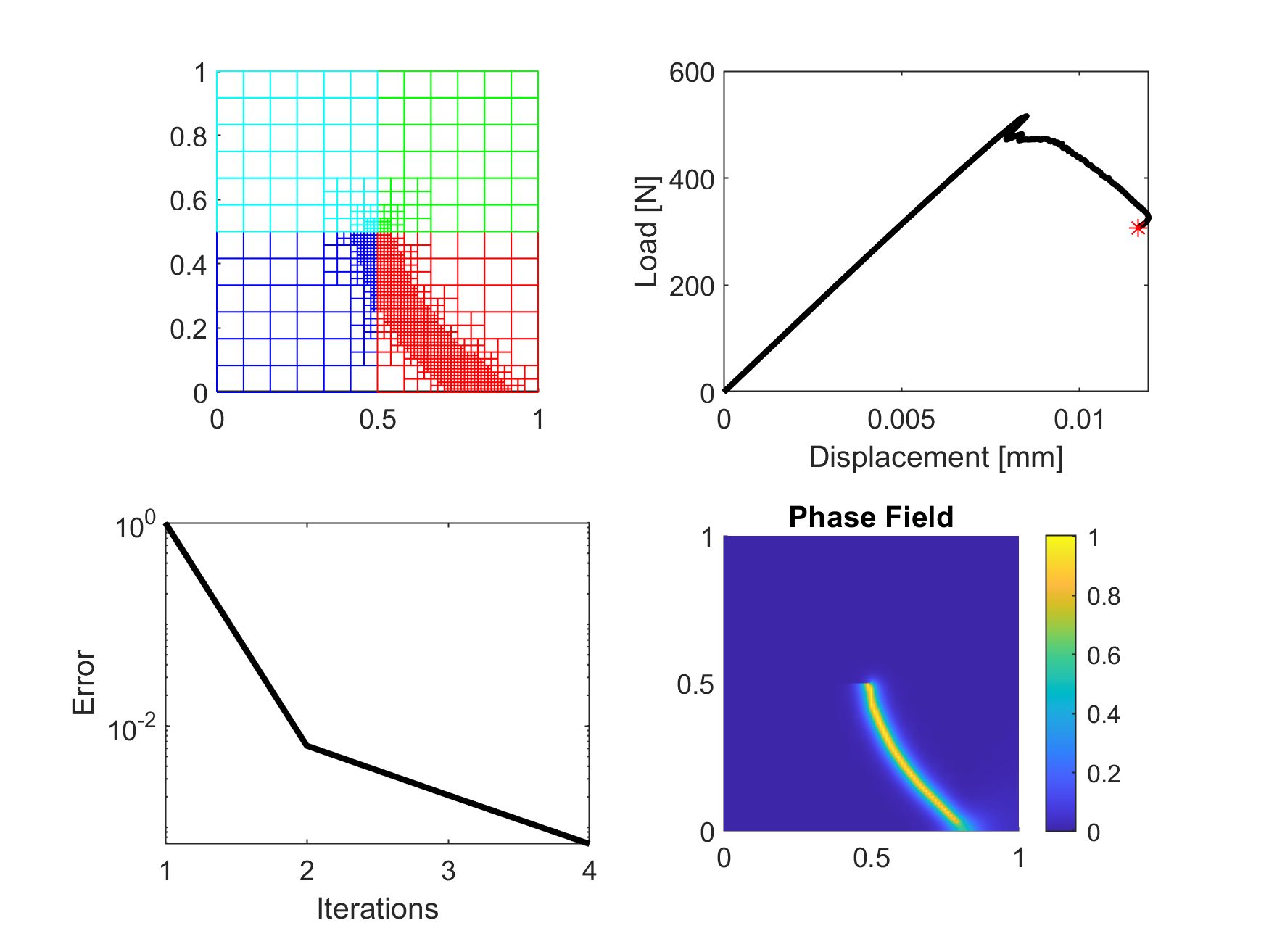}};
    \end{tikzpicture}
    \caption{ }
  \end{subfigure}
  \begin{subfigure}[t]{0.08\textwidth}
    \begin{tikzpicture}
    \begin{axis}[%
    hide axis,
    scale only axis,
    height=1.75\linewidth,
    width=0.01\linewidth,
    point meta min=0.0,
    point meta max=1.0,
    % colormap/bluered,                     % Colormap preset
    colorbar right,                  % Active colorbar
    colorbar sampled,                     % Steps in colorbar
    colorbar style={
        separate axis lines,
        samples=256,                        % Number of steps
    },
    every colorbar/.append style={
    height=\pgfkeysvalueof{/pgfplots/parent axis height},
    ytick={0,0.5,1},yticklabels={0,0.5\:\,$\pf$,1}
    }
    ]
    \addplot [draw=none] coordinates {(0,0)};
    \end{axis}
    \end{tikzpicture}
    \end{subfigure}
  % Mesh refinement plots    
  \begin{subfigure}[t]{0.23\textwidth}
  \centering
    \begin{tikzpicture}[scale=0.5]
    \node[inner sep=0pt] () at (0,0)
    {\includegraphics[width=3.0cm,trim=2.45cm 6.40cm 8.15cm 0.5cm, clip]{Figures/SENS/Ref/Step21.png}};
    \end{tikzpicture}
    \caption{ }
  \end{subfigure}
  %\hfill
  % 2nd
  \begin{subfigure}[t]{0.22\textwidth}
  \centering
    \begin{tikzpicture}[scale=0.5]
    \node[inner sep=0pt] () at (0,0)
    {\includegraphics[width=3.0cm,trim=2.45cm 6.40cm 8.15cm 0.5cm, clip]{Figures/SENS/Ref/Step40.png}};
    \end{tikzpicture}
    \caption{}
  \end{subfigure}
  \begin{subfigure}[t]{0.22\textwidth}
  \centering
    \begin{tikzpicture}[scale=0.5]
    \node[inner sep=0pt] () at (0,0)
    {\includegraphics[width=3.0cm,trim=2.45cm 6.40cm 8.15cm 0.5cm, clip]{Figures/SENS/Ref/Step85.png}};
    \end{tikzpicture}
    \caption{ }
  \end{subfigure}
  \begin{subfigure}[t]{0.22\textwidth}
  \centering
    \begin{tikzpicture}[scale=0.5]
    \node[inner sep=0pt] () at (0,0)
    {\includegraphics[width=3.0cm,trim=2.45cm 6.40cm 8.15cm 0.5cm, clip]{Figures/SENS/Ref/Step126.png}};
    \end{tikzpicture}
    \caption{ }
  \end{subfigure}
  \caption{Figures (a-d) shows the evolution of the phase-field, and the corresponding refined meshes are shown in Figures (e-h) for the single edge notched specimen under shear. The colors in the latter figures represent the different patches in the IGA context.}
  \label{sec5:fig:adap_mesh_SENS}
\end{figure}

Figure \ref{sec5:fig:SENS_lodi} presents the load-displacement curves for the single edge notched specimen under shear.  The different curves correspond to the choice of the degradation function, quadratic and cubic. It is observed that the specimen exhibits a linear pre-peak behaviour with the cubic degradation function for $s<1$. This linear stage is not exhibited by the quadratic degradation function and the cubic degradation function with $s=1$. However, irrespective of the choice of the degradation function, two snap-back behaviours are observed. The first one occurring at the onset of the localization whereas the second snap-back behaviour appears when the crack has reached the bottom edge. Next, in Figure \ref{sec5:fig:SENS_failure}, the phase-field fracture topology at the final step of the analysis is presented. The fracture topology differs from that presented in \cite{Miehe2010}, and the reason lies in the choice of fracture driving $\Psi^f$. In \cite{Miehe2010}, the fracture is driven by the tensile strain energy obtained through spectral decomposition, whereas in this work, the Rankine criterion \cite{wu2017} is adopted. The refined meshes corresponding to the different stages in the evolution of the phase-field is shown in Figure \ref{sec5:fig:adap_mesh_SENS}. 

\newpage

%--------------- Concluding Remarks  -------------------------------%
\section{Concluding Remarks}\label{sec:6}

In this work, we have proposed a robust monolithic solver to accurately forecast the material behavior, which is essential for predicting damage progression and estimating possible failure paths. The literature on variational phase-field based fracture modeling still lacks a reliable, efficient, and simple monolithic solver capable of capturing the pre- and post-peak behaviours accurately. The proposed fully monolithic solver adopts a fracture energy-based arc-length method and an adaptive under-relaxation technique to bridge this gap. The proposed solver utilizes an adaptive mesh refinement scheme using polynomial splines over hierarchical T-meshes (PHT-splines) within the framework of IGA. The PHT-splines possess a very efficient and easy to implement local refinement algorithm, which makes it a right choice for capturing quantities of local interest. The combination of the proposed solver with an adaptive mesh refinement technique could facilitate the application of this approach to more complex structures and with sophisticated constitutive laws. Through the four test cases presented in this work, the crack is allowed to nucleate on its own, and also the solver captures the post-peak snap back effects which is not possible with the alternate minimization solver\cite{Bourdin2007,miehe2010b} and the quasi-monolithic scheme \cite{Heister2015}.

%Further extensions of the present research might include additional verification with more examples involving complex geometries, the adoption of a unified phase field model for both AT1 and AT2 fracture analysis as well as a direct extension to dynamic loading and ductile fracture analysis, the latter being the primary motivation of this work.

Further extensions of this work may include complex multiphysics problems (e.g., porous media, corrosion), dynamic fracture, and the unified phase-field fracture model \cite{wu2017} for quasi brittle fracture. Also, plasticity models could be incorporated for ductile fracture. 

\paragraph{Declaration of Competing Interest:} The authors declare that they have no known competing financial interests or personal relationships that could have appeared to influence the work reported in this paper.

\paragraph{Acknowledgment:} The first author (R.B) is thankful to Elias B\"orjesson for many helpful discussion on the arc-length method. The financial support from the Swedish Research Council for Sustainable Development (FORMAS) under Grant 2018-01249 and the Swedish Research Council (VR) under Grant 2017-05192 is gratefully acknowledged.

\bibliographystyle{unsrt}
\bibliography{ref}

\newpage
\appendix

\section{Solution algorithms}\label{appA:algorithms}

\begin{algorithm}[ht!]
\DontPrintSemicolon
\KwData{Initialisation}
$method \gets `displacement$'\;
$step \gets 0$\;
$\lambda \gets 0$\;
$\Delta\lambda \gets `stepsize$'\;
\ForAll{timesteps}
  {
    $step += 1$\;
    $\lambda += \Delta\lambda$\;
    \While{not Converged}{
    \Switch{method}{
    \Case{`displacement'}
    {Solve Discrete Problem \ref{Problem3}}
    \Case{`Arc-length'}
    {Solve Discrete Problem \ref{Problem4}}
    }
    }
    Compute energy dissipated ($\Delta G$) in the step using (\ref{eqn:sec2:EFracIncr})\;
    \uIf{$\Delta G \geq switchEnergy \:\&\&\: method == `displacement'$}
    {
      $\Delta\lambda = 0$\;
      $\Delta\tau = switchEnergy$\;
    }
    \Else
    {
       \uIf{$iter < optiter$}
       {
         $\Delta\tau *= \Delta\tau$\;
         $\Delta\tau =  min(\Delta\tau,\Delta\tau_{max})$\;
       }
       \Else
       {
         $\Delta\tau = 0.5^{(0.25*(iter-optiter))}$\;
         $\Delta\tau =  min(\Delta\tau,\Delta\tau_{max})$\;
       }
    }
  }\label{endfor}
\caption{Solution scheme without under-relaxation}\label{alg:dispToArc}
\end{algorithm}

\newpage

\begin{algorithm}[ht!]
\DontPrintSemicolon
\KwData{Initialisation}
$method \gets `displacement$'\;
$step \gets 0$\;
$\lambda \gets 0$\;
$\Delta\lambda \gets `stepsize$'\;
\ForAll{timesteps}
  {
    $step += 1$\;
    $\lambda += \Delta\lambda$\;
    \While{not Converged}{
    \Switch{method}{
    \Case{`displacement'}
    {
    Solve Discrete Problem \ref{Problem3}\;
    Update solution fields with relaxation (\ref{eq:underRelaxSolUpdate})
    }
    \Case{`Arc-length'}
    {
    Solve Discrete Problem \ref{Problem4}\;
    Update solution fields with relaxation (\ref{eq:underRelaxSolUpdate})
    }
    }
    }
    \uIf{failed to Converge}
    {
    Revert to old step solution and $\lambda$\;
      \uIf{method == `displacement'}
       {
         $\Delta\lambda /= 10$
       }
       \Else
       {
         $fail += 1$\;
         \uIf{$fail > 2$}
         {
           $\Delta\tau \,/= 2$\;
           $\beta \,/=  1.25$\;
         }
         $\beta \,/=  1.25$\;
      }
    \textbf{break}\;
    }
    $fail \gets 0$\;
    $\beta \,*= 1.05$\;
    Compute energy dissipated ($\Delta G$) in the step using (\ref{eqn:sec2:EFracIncr})\;
    \uIf{$\Delta G \geq switchEnergy \:\&\&\: method == `displacement'$}
    {
      $\Delta\lambda = 0$\;
      $\Delta\tau = switchEnergy$\;
    }
    \Else
    {
       \uIf{$iter < optiter$}
       {
         $\Delta\tau *= \Delta\tau$\;
         $\Delta\tau =  min(\Delta\tau,\Delta\tau_{max})$\;
       }
       \Else
       {
         $\Delta\tau = 0.5^{(0.25*(iter-optiter))}$\;
         $\Delta\tau =  min(\Delta\tau,\Delta\tau_{max})$\;
       }
    }
  }
\caption{Solution scheme with under-relaxation}\label{alg:dispToArcUnderRelax}
\end{algorithm}

\section{Sensitivity w.r.t maximum dissipation allowed in a step and refinement indicator}\label{appA:sensitivity}

In this section, the sensitivity of the load-displacement curves w.r.t the maximum dissipated energy $\Delta\tau_{max}$ [N] and refinement indicator ($\pf_{th.}$) is investigated. To this end, the single edge notched specimen under tension (SENT) from Section \ref{sec5:SENtension} is considered. All aspects of the model (geometry, material parameters and loading conditions) remain same, the only variation being in $\Delta \tau_{max}$ and the degradation function set to quadratic. 

\begin{figure}[!ht]
\begin{subfigure}[t]{0.45\textwidth}
\centering
\begin{tikzpicture}[thick,scale=0.95, every node/.style={scale=0.9}]
    \begin{axis}[legend style={draw=none}, legend columns = 2,
      transpose legend, ylabel={Load\:[N]},xlabel={Displacement\:[mm]}, xmin=0, ymin=0, xmax=0.0065, ymax=1000, yticklabel style={/pgf/number format/.cd,fixed,precision=2},
                 every axis plot/.append style={very thick}]
    \pgfplotstableread[col sep = comma]{./Data/VaryDis/lodi_dTau0_00625.txt}\Adata;
    \pgfplotstableread[col sep = comma]{./Data/VaryDis/lodi_dTau0_0125.txt}\Bdata;
    \pgfplotstableread[col sep = comma]{./Data/SENT/lodi_lby2pf0_2.txt}\Cdata;
    \pgfplotstableread[col sep = comma]{./Data/VaryDis/lodi_dTau0_05.txt}\Ddata;
    \addplot [ 
           color=black, 
%           only marks, 
%           mark=*, 
%           mark size=0.25pt, 
         ]
         table
         [
           x expr=\thisrowno{1}, 
           y expr=\thisrowno{0}
         ] {\Adata};
         \addlegendentry{$\Delta\tau_{max} = 0.00625$}
    \addplot [ 
           color=red, 
%           only marks, 
%           mark=*, 
%           mark size=0.25pt, 
         ]
         table
         [
           x expr=\thisrowno{1}, 
           y expr=\thisrowno{0}
         ] {\Bdata};
         \addlegendentry{$\Delta\tau_{max} = 0.0125$}     
    \addplot [ 
           color=blue, 
%           only marks, 
%           mark=*, 
%           mark size=0.25pt, 
         ]
         table
         [
           x expr=\thisrowno{1}, 
           y expr=\thisrowno{0}
         ] {\Cdata};
         \addlegendentry{$\Delta\tau_{max} = 0.025$}
    \addplot [ 
           color=green, 
%           only marks, 
%           mark=*, 
%           mark size=0.25pt, 
         ]
         table
         [
           x expr=\thisrowno{1}, 
           y expr=\thisrowno{0}
         ] {\Ddata};
         \addlegendentry{$\Delta\tau_{max} = 0.05$}       
    \end{axis}
    \end{tikzpicture}
    \caption{ }
    \label{appA:fig:SENT_dtaumax}
\end{subfigure}
\hfill
  \begin{subfigure}[t]{0.45\textwidth}
  \centering
    \begin{tikzpicture}[thick,scale=0.95, every node/.style={scale=0.9}]
    \begin{axis}[legend style={draw=none}, legend columns = 2,
      transpose legend, ylabel={Load\:[N]},xlabel={Displacement\:[mm]}, xmin=0, ymin=0, xmax=0.0065, ymax=1000, yticklabel style={/pgf/number format/.cd,fixed,precision=2},
                 every axis plot/.append style={very thick}]
    \pgfplotstableread[col sep = comma]{./Data/SENT/lodi_lby2pf0_1.txt}\Adata;
    \pgfplotstableread[col sep = comma]{./Data/SENT/lodi_lby2pf0_2.txt}\Bdata;
    \pgfplotstableread[col sep = comma]{./Data/SENT/lodi_lby2pf0_4.txt}\Cdata;
    \pgfplotstableread[col sep = comma]{./Data/SENT/lodi_FixMesh.txt}\Ddata;
    \addplot [ 
           color=black, 
%           only marks, 
           mark=*, 
           mark size=0.25pt, 
         ]
         table
         [
           x expr=\thisrowno{1}, 
           y expr=\thisrowno{0}
         ] {\Adata};
         \addlegendentry{$\pf_{th.} = 0.1$}
    \addplot [ 
           color=red, 
%           only marks, 
           mark=*, 
           mark size=0.25pt, 
         ]
         table
         [
           x expr=\thisrowno{1}, 
           y expr=\thisrowno{0}
         ] {\Bdata};
         \addlegendentry{$\pf_{th.} = 0.2$}     
    \addplot [ 
           color=blue, 
%           only marks, 
           mark=*, 
           mark size=0.25pt, 
         ]
         table
         [
           x expr=\thisrowno{1}, 
           y expr=\thisrowno{0}
         ] {\Cdata};
         \addlegendentry{$\pf_{th.} = 0.4$}     
    \addplot [ 
           color=green, 
%           only marks, 
           mark=*, 
           mark size=0.25pt, 
         ]
         table
         [
           x expr=\thisrowno{1}, 
           y expr=\thisrowno{0}
         ] {\Ddata};
         \addlegendentry{Fixed Mesh}          
    \end{axis}
    \end{tikzpicture}
    \caption{ }
    \label{appA:fig:SENT_ref}
  \end{subfigure}
  \caption{Figure (a) presents the load-displacement curves for the single edge notched specimen under tension. The legend entries correspond to the choice of prescribed maximum dissipation energy $\Delta\tau_{max}$. Figure (b) presents the load-displacement curves for the single edge notched specimen under tension. The legend entries correspond to the choice of solution techniques.}
\end{figure}
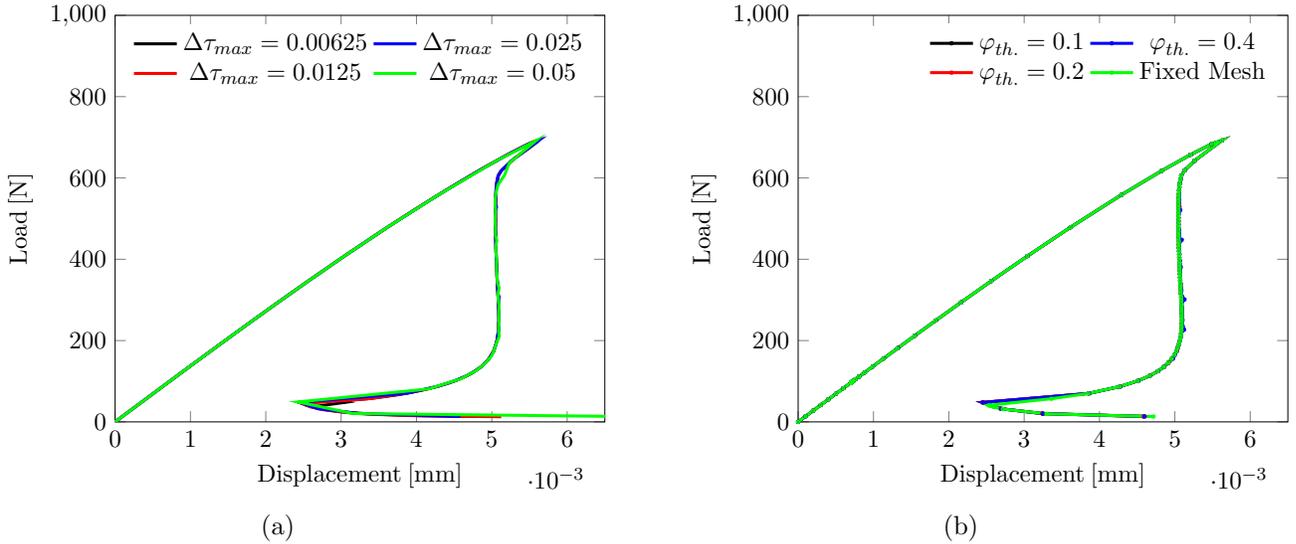

Figure \ref{appA:fig:SENT_dtaumax} presents the load-displacement curves obtained using the proposed solver for different values of $\Delta\tau_{max}$. For the values chosen in this study, the load-displacement curves are similar. This indicates that the relevant features (snap-back behaviour) could be sufficiently captured even with a larger dissipation steps, for instance $\Delta\tau_{max} = 0.05$ [N]. 

The sensitivity of the load-displacement curves obtained using the proposed solver with different refinement indicator is studied with three different threshold phase-field values. An element is marked for refinement once this threshold value is exceeded. Figure \ref{appA:fig:SENT_ref} presents the load-displacement curves obtained for the different phase-field threshold values along with that obtained using a fixed mesh. It is observed that the curves are similar. This indicates that the proposed monolithic solver does not exhibit any bias w.r.t. adaptive mesh refinement techniques.

%\begin{table}[ht!]
%    \centering
%    \begin{tabular}{llllll} \hline
%  $\Delta\tau$ & Total Steps & Solution steps & AMR steps & Fail steps & Avg. Iter.  \\ \hline
%  0.00625 & 313 & 286 & 22 & 5 & 3.131 \\
%  0.0125 & 193 & 164 & 22 & 7 & 4.3472 \\
%  0.025 & 170 & 101 & 58 & 11 & 7.8588 \\ 
%  0.05 & 101 & 70 & 22 & 9 & 13.7525 \\ \hline
%  \end{tabular}
%    \caption{Performance of the monolithic solver w.r.t varying dissipation allowed in a single step. Solution steps refer to those steps which resulted in a converged solution. AMR steps also result in convergence, however, the step is re-computed on the new mesh. Fail steps refer to those steps where the solver failed to achieve convergence. The average iterations (Avg. Iter.) is computed as total iterations by total steps. }
%    \label{appA:table:SENT_varydTau_performance}
%\end{table}

\section{Comparison with Alternative minimization solver and Quasi-Newton Raphson method \cite{Heister2015}}\label{appA:comparison}

In this section, a comparison of the proposed monolithic solver with alternative minimization solver \cite{Bourdin2007,Miehe2010,miehe2010b} and the quasi Newton-Raphson method proposed by \cite{Heister2015} is carried out. The comparison is based only on the load-displacement curves obtained from the respective solvers. The alternate minimization solver solves the displacement and the phase-field sub-problems alternatively until a certain tolerance criterion is met. The quasi Newton Raphson approach proposed in \cite{Heister2015} adopts a linear extrapolation of the phase-field variable for the momentum balance equation. Although, the extrapolation strategy is questionable, it yields a convex energy functional, resulting in a better convergence behaviour of the Newton Raphson method.

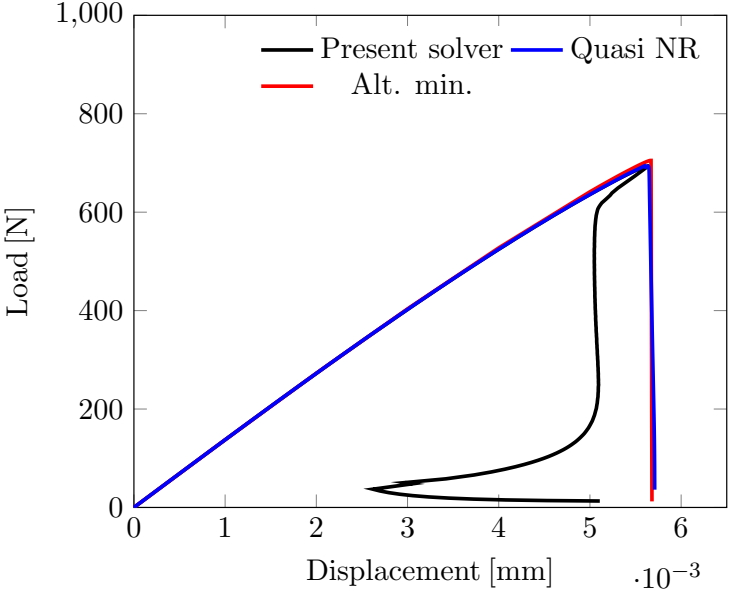
\begin{figure}[!ht]
 \centering
    \begin{tikzpicture}[thick,scale=1.15, every node/.style={scale=0.9}]
    \begin{axis}[legend style={draw=none}, legend columns = 2,
      transpose legend, ylabel={Load\:[N]},xlabel={Displacement\:[mm]}, xmin=0, ymin=0, xmax=0.0065, ymax=1000, yticklabel style={/pgf/number format/.cd,fixed,precision=2},
                 every axis plot/.append style={very thick}]
    \pgfplotstableread[col sep = comma]{./Data/VaryDis/lodi_dTau0_00625.txt}\Adata;
    \pgfplotstableread[col sep = comma]{./Data/CompAMQNR/lodiAM.txt}\Bdata;
    \pgfplotstableread[col sep = comma]{./Data/CompAMQNR/lodiQNR.txt}\Cdata;
    \addplot [ 
           color=black, 
%           only marks, 
%           mark=*, 
%           mark size=0.25pt, 
         ]
         table
         [
           x expr=\thisrowno{1}, 
           y expr=\thisrowno{0}
         ] {\Adata};
         \addlegendentry{Present solver}
    \addplot [ 
           color=red, 
%           only marks, 
%           mark=*, 
%           mark size=0.25pt, 
         ]
         table
         [
           x expr=\thisrowno{1}, 
           y expr=\thisrowno{0}
         ] {\Bdata};
         \addlegendentry{Alt. min.}     
    \addplot [ 
           color=blue, 
%           only marks, 
%           mark=*, 
%           mark size=0.25pt, 
         ]
         table
         [
           x expr=\thisrowno{1}, 
           y expr=\thisrowno{0}
         ] {\Cdata};
         \addlegendentry{Quasi NR}
    \end{axis}
    \end{tikzpicture}
\caption{Figure presents the load-displacement curves for the single edge notched specimen under tension. The legend entries correspond to the choice of solution techniques.}
\label{appB:fig:SENT_lodi_compare}
\end{figure}

For the comparative study, the SENT model from Section \ref{sec5:SENtension} is chosen. For both the alternate minimization solver and the quasi Newton Raphson method, the fracture irreversibility is enforced using the penalization method \cite{Gerasimov2016,GERASIMOV2019990}. Figure \ref{appB:fig:SENT_lodi_compare} presents the load-displacement curves obtained using the solver proposed in this work, and the alternate minimization solver and the quasi-Newton Raphson method. It is observed that all methods/solvers predict a similar peak load and pre-peak behaviour. However, the post-peak behaviours are different. The proposed monolithic solver equipped with the arc-length method is able to predict snap-back behaviour, which is not possible with conventional alternate minimization solver and the quasi-Newton Raphson method.

\newpage

\end{document}